\documentclass[10pt,twocolumn,twoside]{IEEEtran}
\usepackage{amsfonts,amsmath,amssymb,amsthm}
\usepackage[linesnumbered, ruled]{algorithm2e}
\usepackage{times}
\usepackage{bm}
\usepackage{balance}
\usepackage{color}
\usepackage{graphicx}
\usepackage{setspace}
\usepackage{comment}
\usepackage{cases}
\graphicspath{{./eps/}}
\usepackage{here}
\usepackage{ascmac} 
\usepackage{array}
\usepackage{float}
\usepackage{cite}
\usepackage{url}
\usepackage{subfigure}
\newtheorem{mydef}{Definition}
\newtheorem{mylem}{Lemma}

\newtheorem{mythm}{Theorem}
\newtheorem{myas}{Assumption}
\newtheorem{myrem}{Remark}

\newcommand{\rfig}[1]{Fig.\,\ref{#1}} 
\newcommand{\rapp}[1]{Appendix~\ref{#1}} 
\newcommand{\req}[1]{\eqref{#1}}

\newcommand{\rlem}[1]{Lemma\,\ref{#1}}

\newcommand{\rrem}[1]{Remark\,\ref{#1}}
\newcommand{\rsec}[1]{Section\,\ref{#1}}

\newcommand{\ras}[1]{Assumption\,\ref{#1}}
\newcommand{\ralg}[1]{Algorithm\,\ref{#1}}

\newcommand{\rline}[1]{line\,\ref{#1}}

\newcommand{\rthm}[1]{Theorem\,\ref{#1}}

\SetKwInOut{Parameter}{parameter}
\newcommand{\qedwhite}{\hfill \ensuremath{\Box}}

\begin{document}
\title{Distributed $\ell_1$-state-and-fault estimation \\ for Multi-agent systems}
\author{\large Kazumune~Hashimoto, Michelle~Chong and Dimos~V.~Dimarogonas,~\IEEEmembership{Senior Member,~IEEE}
\thanks{The authors are with the School of Electrical Engineering and Computer Science, KTH Royal Institute of Technology, 10044 Stockholm, Sweden. This work is supported by the Knut and Alice Wallenberg Foundation.}
}
\maketitle

\begin{abstract}
In this paper, we propose a distributed state-and-fault estimation scheme for multi-agent systems. 
The estimator is based on an $\ell_1$-norm optimization problem, which is inspired by sparse signal recovery in the field of compressive sampling. Two theoretical results are given to analyze the correctness of our approach. 
First, we provide a necessary and sufficient condition such that the state and fault signals are correctly estimated. The result presents a fundamental limitation of the algorithm, which shows how many faulty nodes are allowed to ensure a correct estimation. Second, we analyze how the estimation error grows over time by showing that the upper bound of the estimation error depends on the previous state estimate and the number of faulty nodes. 
An illustrative example is given to validate the effectiveness of the proposed approach. 
\end{abstract}

\section{Introduction}
The problem of controlling multi-agent systems arises in many practical applications, such as vehicle platooning, formation control of autonomous vehicles and satellites, cooperative control of robots, power networks, to name a few \cite{magnus}. In multi-agent systems, designing an automatic framework of detecting \textit{anomaly} or \textit{faults} has attracted much attention in recent years, see, e.g., \cite{faultsurvey}. The motivation is two fold. First, since a multi-agent system consists of many interactions among subsystems, only a few anomalous behaviors will provide a significant impact on the overall system's performance. For example, in vehicle platooning, anomalous behaviors such as sudden braking or acceleration of a leading vehicle will affect behaviors of the following vehicles, which may cause traffic accidents or traffic jam. Second, due to the existence of many interactions over the communication network, a multi-agent system may be vulnerable to \textit{cyber attacks}. For example, actuators of some nodes may be hacked by attackers over communication network, so that they can operate the controllers to enforce some nodes to be anomalous. Therefore, detecting and estimating faults in multi-agent system at an early stage is necessary to maintain safe and reliable real-time operations. 

So far, various formulations and approaches have been proposed to detect and estimate faults in multi-agent systems, see, e.g., \cite{ding,fdi5,fdi0,fdi1,fdi2,fault_detect1,fdi3,fdi4,fdi7,fdi9,fdi10}. 
Early works focus on designing detection schemes in a \textit{centralized} manner, in which all nodes should communicate with a central unit for detecting faulty signals, see, e.g., \cite{ding,fdi5}. In recent years, there has been a growing attention in designing \textit{decentralized} or \textit{distributed} fault detection and estimation schemes, see, e.g., \cite{fdi0,fdi1,fdi2,fault_detect1,fdi3,fdi4,fdi9,fdi7,fdi10}. 
For example, \cite{fdi0} and \cite{fdi1} proposed Unknown Input Observers (UIO) for detecting faults in a decentralized and distributed manner, respectively. In \cite{fdi2}, a distributed fault estimation scheme was proposed by using $H_\infty$ optimization, which involves sensitivity of faults and robustness of disturbances. In \cite{fdi10}, a consensus-based decentralized observer is designed for fault detection. 
Moreover, a distributed fault detection scheme based on UIO for multi-agent Linear Parameter Varying (LPV) systems is proposed in \cite{fdi9}. 

In this paper, we are interested in developing a \textit{distributed} fault estimation strategy for multi-agent systems. In particular, we consider the situation where each node is able to measure \textit{relative state information} with respect to its neighbors. 
Some approaches have been already proposed for this problem set-up, see, e.g., \cite{fault_detect1,fdi3,fdi4,fdi7}, where the references include centralized, decentralized, and distributed schemes. For example, in \cite{fdi4} a coordinate transformation is introduced to extract the observable subspace, and a sliding mode decentralized observer is designed to estimate faulty signals. In \cite{fdi7}, a distributed fault detection scheme is proposed by employing Linear Matrix Inequalities (LMIs) on the sensitivity of residuals and robustness against system noise. In \cite{fdi3}, a fault detection scheme is developed by introducing a set-valued observer, which is applied to two subsystems that are decomposed by left-coprime factorization. In addition, a UIO-based event-triggered communication protocol was proposed in \cite{fault_detect1}, where each node monitors the neighboring nodes if they are subject to faults. 

The contribution of this paper is to present a novel distributed fault estimation scheme for multi-agent systems based on relative state measurements. 
The estimator is different from previous works in \cite{fault_detect1,fdi3,fdi4,fdi7}, and draws inspirations from compressive sampling \cite{compressivesampling}. The estimator can handle the case when multiple faults arise. Thus, the approach is advantageous over \cite{fault_detect1,fdi3,fdi4} where neither systematic estimation scheme for multiple faults nor theoretical analysis are given. Moreover, we provide a quantitative analysis of how many faulty nodes can be tolerated to provide a correct estimation. While \cite{fdi7} gives a sufficient condition for designing a suitable estimator for multiple faults, our result may be more specific and intuitive, giving us a necessary and sufficient condition when the estimator provides correct estimations. 

Specifically, the contribution of this paper is as follows: 



\begin{enumerate}
\item We {provide} a fault-and-state estimator for multi-agent systems based on an \textit{$\ell_1$-norm optimization problem}, which is inspired by sparse signal recovery in the field of compressive sampling \cite{compressivesampling}. The optimization problem is formulated as a Basis Pursuit \cite{basispursuit2}, in which several numerical solvers can be employed to solve the problem efficiently. 

\item We provide two quantitative analysis for the correctness of the estimator. 
First, we provide a necessary and sufficient condition such that states and faulty signals are correctly estimated for all times. Indeed, this illustrates a fundamental limitation of the algorithm, which shows how many faulty nodes can be tolerated to provide correct estimation. 
Second, we consider the case when the estimation error is present for some time due to, e.g., numerical errors of solving the optimization problem, and provide error bounds for the state-and-fault estimation. 
In particular, we show that the estimation error is bounded depending on not only the previous state estimate but also the number of faulty nodes. The result also shows how the estimation error of states grows over time.

\item We provide a distributed implementation to solve the $\ell_1$-norm optimization problem. 
To this aim, we employ the approach presented in \cite{basispursuit}, which utilizes the alternating direction method of multipliers (ADMM). As with \cite{basispursuit}, the equivalence between the centralized and the distributed solution is shown. 

\end{enumerate}

{Since the approach employs the $\ell_1$-norm optimization, the analysis provided in this paper may be relevant to secure estimation for Cyber Physical Systems (CPSs), see, e.g., \cite{fawzi,shoukry2016,michelle}. 
The problem formulation presented here differs from these previous works in the following sense. While previous works consider that both the initial state of the system and attack (fault) signals are unknown and are to be estimated, in this paper we consider that attack (fault) signals are unknown, but the initial state of the system is assumed to be \textit{known}. 
The assumption that the initial state is known limits the applicability of our approach. However, examples of physical systems where the initial state is known include multi-vehicle systems, in which the initial position of each vehicle can be known by utilizing the GPS, such as the case we consider in \rsec{simulation_example}. 
More importantly, while the initial state is assumed to be known, we deal with the case when the dynamics of the plant is given by an \textit{unobservable} system, 
which is due to that only relative state measurements are available. 
Note that the analysis for unobservable systems where the initial state is known has not been investigated in the previous works of \cite{fawzi,shoukry2016,michelle}. 
An interesting observation may be that, while our problem formulation and the previous ones are different as above, the corresponding results are relevant, since they both provide a connection between the number of faulty nodes (attacked sensors) and the correctness of the estimation.

Our \textit{centralized} estimator is analogous to the approach in \cite{fdi11}. 
However, our result is novel with respect to \cite{fdi11} in the following sense. 
First, we provide a more detailed analysis for multi-agent systems with relative state measurements, which derives a concrete connection between the number of faulty nodes and the correctness of the estimation. 
As we will see later, this is achieved by analyzing the null-space property of the measurement matrix. Second, we formulate a \textit{distributed} implementation to solve the $\ell_1$-norm optimization problem, in which each node is able to estimate states and faults for all nodes by coordinating only with its neighbors.

The remainder of this paper is organized as follows. 
In Section~II, we present some preliminaries on graph theory and the $\ell_1$-norm optimization problem. In Section~III, the problem formulation and approach are given. In Section~IV, we provide analysis of the estimator. In Section~V, we provide analysis of the estimator in the presence of estimation errors at previous time. In Section~VI, a distributed implementation to solve the estimation problem is given. In Section~VII, an illustrative example shows the effectiveness of the proposed approach. Conclusions are given in Section~VIII. \\

\noindent

\noindent
\textbf{Notation:} Let $\mathbb{R}$, $\mathbb{R}_+$, $\mathbb{N}$, $\mathbb{N}_+$ be the set of reals, {positive reals, non-negative integers}, and {positive integers}, respectively.  
Given $x\in\mathbb{R}^n$, denote by $x^{(i)}$ the $i$-th component of $x$. Given $x\in\mathbb{R}^n$, let $\|x\|$ be the Euclidean norm of $x$. 
Moreover, let $\|x\|_1$ be 
the $\ell_1$-norm of $x$, i.e., $\|x\|_1 = |x^{(1)}| +|x^{(2)}| + \cdots + |x^{(n)}|$. Given a matrix $A \in\mathbb{R}^{n\times m}$, we denote by $A^{(i, j)}$ the $(i,j)$-component of $A$. 
Given $M\in\mathbb{N}_+$, denote by $1_M \in \mathbb{R}^M$ a $M$-dimensional vector whose components are $1$. 
Given a set $T$, denote by $|T|$ the cardinality of $T$. 
For given $T\subseteq \{1, \ldots, n\}$ and $x\in\mathbb{R}^n$, $x$ is called $T$-sparse if all components indexed by the complement of $T$ is $0$, i.e., $x^{(i)} = 0$, $\forall i\in\{1, \ldots, n\}\backslash T$. 
Given $x\in\mathbb{R}^n$ and $T\subset \{1, \ldots, n\}$, denote by $x_T \in \mathbb{R}^{|T|}$ a vector from $x$ by extracting all components indexed by $T$. 
Given $A \in\mathbb{R}^{m\times n}$, denote by ${\rm rank}(A)$ the {rank} of $A$. Moreover, denote by ${\rm ker}(A)$ the {null space} of $A$, i.e., 
${\rm ker}(A) = \{ x \in \mathbb{R}^n : A x = 0 \}$.

\section{Preliminaries}
\subsection{Graph theory}\label{pre_graphsec}
Let ${\cal G} = ({\cal V}, {\cal E})$ denote a directed graph, where ${\cal V} = \{1, 2, \ldots, M\}$ is the set of nodes and ${\cal E} = \{e_1, e_2, \ldots, e_N\} \subseteq {\cal V} \times {\cal V}$ is the set of edges. The graph is strongly connected if for every pair of nodes there exists a path between them. The graph is weakly connected if for every pair of nodes, when removing all orientations in the graph, there exists a path between them. For a given $i\in{\cal V}$, let ${\cal N}_i \subset {\cal V}$ be the set of neighboring nodes of $i$, i.e., ${\cal N}_i = \{j\in {\cal V} : e = \{i, j\} \in {\cal E}\}$. The incidence matrix $D = D(G)$ is the $\{0, \pm 1\}$ matrix, where $D^{(i,j)} = 1$ if node $i$ is the head of the edge $e_j$, $D^{(i,j)} = -1$ if node $i$ is the tail of the edge $e_j$, and $D^{(i,j)} = 0$ otherwise. 
For a given ${\cal G} = ({\cal V}, {\cal E})$ with $M$ number of nodes, it is shown that ${\rm rank} (D) = M-1$ if ${\cal G}$ is weakly connected. The null-space of the incidence matrix is given by ${\rm ker}(D^\mathsf{T}) = \gamma {1}_M$, where $\gamma \in \mathbb{R}$ (see, e.g., Chapter~2 in \cite{graph}). 

\subsection{The Null-Space Property and $\ell_1$-norm optimization} 
In what follows, we provide the notion of {Null-Space Property} (NSP) and its useful result for the $\ell_1$-norm optimization problem. 
\begin{mydef}[The Null-Space Property]
\normalfont
For given $A \in \mathbb{R}^{m\times n}$ and $T \subseteq \{1, 2, \ldots, n\}$, $A$ is said to satisfy the Null-Space Property (NSP) for $T$ (or $T$-NSP for short), if for every $v \in {\rm ker} (A) \backslash \{0\}$, it holds that 
\begin{equation}
\|v_T\|_1 < \|v_{T^c}\|_1, 
\end{equation}
where $T^c = \{1, \ldots, n\}\backslash T$. 
 \qedwhite
\end{mydef}
The NSP is a key property to check whether the sparse signal can be reconstructed based on the $\ell_1$-norm optimization problem: 


\begin{mythm}[$\ell_1$-reconstruction theorem]\label{reconstruct}
\normalfont
For given $A \in \mathbb{R}^{m\times n}$ and $T \subseteq \{1, \ldots, n\}$, every $T$-sparse vector $x_0 \in \mathbb{R}^n$ is a unique solution to the following optimization problem: 
\begin{equation}\label{l1_original}
\underset{x}{{\min}}\ \| x \| _1 \ \ \ {\rm s.t.,}\ \ Ax_0 = A x, 
\end{equation}
if and only if $A$ satisfies $T$-NSP. \qedwhite
\end{mythm}
The above theorem indicates that every $T$-sparse vector can be reconstructed by solving the $\ell_1$-norm optimization problem in \req{l1_original}, if and only if the matrix $A$ satisfies $T$-NSP. The proof of \rthm{reconstruct} follows the same line as \cite{sparse1} and is given in the Appendix. 

\section{Problem formulation}\label{problem_formulation_sec}
\subsection{Dynamics} \label{dynamics_sec}
We consider a network of $M$ inter-connected nodes, 
which is modeled by a graph ${\cal G} = ({\cal V}, {\cal E})$, where ${\cal V} =\{1 , 2, \ldots, M\}$ is the set of nodes and ${\cal E} \subseteq {\cal V} \times {\cal V}$ is the set of edges. Here, each edge $e = \{i, j\} \in {\cal E}$ indicates the sensing and communication capabilities of node $i$ with respect to $j$. More specifically, if $\{ i, j \} \in {\cal V}$ node $i$ is able to measure relative state information with respect to $j$, as well as to transmit it to $j$. 
We assume that the network is weakly connected. Moreover, we assume that there exists one leader and $M-1$ followers; without loss of generality, the node labeled with $1\in {\cal V}$ indicates the leader and the others are the followers. 
The dynamics for node $i \in {\cal V}$ is described by the following discrete-time system: 
\begin{align}\label{state_equation}
x_i (k) & = A_i x_i(k-1) + B_i u_i (k-1) + f_i (k), 
\end{align}
for $k\in\mathbb{N}$, where $x_i \in\mathbb{R}^n$ is the state, $u_i \in \mathbb{R}^m$ is the control input, and $f_i \in\mathbb{R}^n$ is the signal indicating the occurrence of faults, i.e., $f_i (k) \neq 0$ when node $i$ is subject to a fault at $k$, and $f_i (k) = 0$ otherwise. {Here, we do not assume any statistical properties of the faulty signals}, and these can be viewed as system faults on the dynamics in \req{state_equation}, or actuator faults caused by physical effects or cyber attacks over the communication network. 
Note that the state and input dimensions are $n$ and $m$ for all nodes. We assume that the control input is known to node $i$ for all $k\in\mathbb{N}$.
Let ${\cal I}_k \subseteq \{ 1, \ldots, M\} \cup \emptyset$ be the \textit{unknown} set of nodes that are subject to faults at $k\in\mathbb{N}$, i.e., 
\begin{equation}\label{anomaly_detection}
{{\cal I}}_k = \left \{i\in\{1, \ldots, M\} : {f}_i (k) \neq 0 \right \}. 
\end{equation}
Note that the set ${\cal I}_k$ can be time-varying, i.e., the set of faulty nodes is possibly changing over time. 

The output equation depends on whether the node is the leader or one of the followers. 
For the followers, each node is able to measure the \textit{relative state information} from its neighbors by using on-board sensors, i.e., for all $i\in {\cal V}\backslash \{1\}$, 
\begin{equation}\label{relative_measurement}
y_{ij} (k) = x_i (k) - x_j (k), \ \ j \in {\cal N}_i
\end{equation}
for $k\in\mathbb{N}$, where $y_{ij} (k)$, $j \in {\cal N}_i$ is the output measurement. 
The set of outputs collected at node $i$ can be expressed as 
\begin{equation}\label{output_follower}
y_{i} (k) = C_i x (k),  
\end{equation}
where $y_i (k) \in \mathbb{R}^{n|{\cal N}_i|}$ is the vector collecting the measurement outputs for all $j \in {\cal N}_i$, $C_i$ is the matrix defined through \req{relative_measurement}, and $x(k) = [x_1(k)^\mathsf{T}\ x_2(k)^\mathsf{T}\ \ldots\  x_M(k)^\mathsf{T}]^\mathsf{T}$. 

Regarding the measurement outputs for the leader, we consider the following two cases. First, not only can the leader measure the relative state information from its neighbors but also its own state information, i.e., 
\begin{align}
y_{1j} (k) &= x_1 (k) - x_j (k), \ \ j \in {\cal N}_1, \label{measurement_leader} \\ 
y_{11} (k) &= x_1 (k), \label{measurement_leader2}
\end{align}
where $y_{1j} (k)$, $j \in {\cal N}_1$ and $y_{11} (k)$ are the output measurements. 
If the leader is able to obtain both \req{measurement_leader} and \req{measurement_leader2}, we say that the leader is in \textit{active mode}. Second, when the leader can only measure the relative state info from its neighbors, i.e., 
\begin{align}\label{measurement_leader3}
y_{1j} (k) &= x_1 (k) - x_j (k), \ \ j \in {\cal N}_1, 
\end{align}
then we say that the leader is in \textit{non-active mode}. In practice, the leader's mode indicated above changes over time depending on physical environments. For example, if the leader is in the outdoor environment that can utilize GPS to locate the state, the measurement in \req{measurement_leader2} is available and thus the leader is in active mode. 
On the other hand, if the leader enters an indoor environment in which it cannot utilize GPS due to signal loss caused by the presence of walls in the building (e.g., inside the tunnel, underground, etc.), the measurement in \req{measurement_leader2} is \textit{not} available and the leader is in non-active mode. 
To indicate whether the leader is active or non-active, let $a_{1} (k) \in \{0, 1\}$ be given by 
\begin{numcases}
{a_1 (k) = }
	0, \ \ {\rm if\ the\ leader\ is\ in\ nonactive \ mode}.\notag \\
	1, \ \ {\rm if\ the\ leader\ is\ in\ active \ mode.} \notag 
\end{numcases}
We further assume the following. 
\begin{myas}\label{initial_assumption}
\normalfont
The leader is in active mode at the initial time, i.e., $a_1(0)=1$. Moreover, there is no fault at $k=0$, i.e., ${\cal I}_0 = \emptyset$. 
\end{myas}
Note that the leader can either be in active or non-active mode for all $k\in\mathbb{N}\backslash\{0\}$. 
As we will see in the analysis of \rsec{analysis_no_error}, 
the assumption of being in an active mode at $k=0$ leads to the correct fault and state estimation for all the time steps afterwards. As with \req{output_follower}, the set of outputs for the leader can be expressed as 
\begin{align}
y_{1} (k) &= C_{10} x (k),\ \ \ \ {\rm if}\ a_1 (k) = 0, \\
y_{1} (k) &= C_{11} x (k),\ \ \ \ {\rm if}\ a_1 (k) = 1,  \label{output_leader}
\end{align}
where  $C_{10}$ is the matrix obtained by collecting all measurements from \req{measurement_leader3}, and 
$C_{11}$ is the matrix obtained from \req{measurement_leader} and \req{measurement_leader2}. 

Based on the above, the overall dynamics for all nodes can be expressed as 
\begin{align}
x(k) &= A x(k-1) + B u(k-1) + f(k), \label{overall_state}
\end{align}
for $k\in\mathbb{N}$, where 
\begin{align}
x(k) &= [x_1(k)^\mathsf{T}\ x_2(k)^\mathsf{T}\ \ldots x_M(k)^\mathsf{T}]^\mathsf{T}, \\
u(k) &= [u_1(k)^\mathsf{T}\ u_2(k)^\mathsf{T}\ \ldots u_M(k)^\mathsf{T}]^\mathsf{T}, \\
f(k) &= [f_1(k)^\mathsf{T}\ f_2(k)^\mathsf{T}\ \ldots f_M(k)^\mathsf{T}]^\mathsf{T}, 
\end{align}
and $A, B$ are the matrices with appropriate dimensions. The measurement output is given by 
\begin{numcases}
{y (k) = }
	C _{0} x(k), \ \ {\rm if}\ a_1 (k) = 0 \label{overall_output2} \\
	C _{1} x(k), \ \ {\rm if\ }a_1 (k) = 1,  \label{overall_output}
\end{numcases}
where $y(k) = [y_1(k)^\mathsf{T}\ y_2(k)^\mathsf{T}\ \ldots y_M(k)^\mathsf{T}]^\mathsf{T}$ and 
\begin{equation}\label{cpartition}
C_0 =  \left [
\begin{array}{c}
 C_{10}  \\
 C_2 \\
 \vdots \\
 C_M
\end{array}
\right ],\ \  C_1 =  \left [
\begin{array}{c}
 C_{11}  \\
 C_2 \\
 \vdots \\
 C_M
\end{array}
\right ]. 
\end{equation}

Note that the pair $(A, C_0)$ is \textit{un-observable} when all nodes are identical (i.e., $A_1 = A_2 = \cdots = A_M$), see e.g., \cite{fdi4,fdi3,fault_detect1}, and such case is included in the above system when the leader is in non-active mode. \\




\noindent
\textit{(Example)} : Suppose that ${\cal V} = \{1, 2, 3\}$ and $\{1, 2\}, \{2, 3 \} \in {\cal E}$, and the dynamics of node $i$ is scalar and given by $x_i (k+1) = x_i (k) + f_i (k)$, $\forall i\in\{1, 2, 3\}$. Then, we have $A = I_3$, $B=0$ and 
\begin{align}\label{c_example}
C_0 = \left [
\begin{array}{ccc}
1  &  -1 & 0 \\
0  &  1  & -1 \\
\end{array}
\right ], \ \ \ 
C_1 = \left [
\begin{array}{ccc}
 1 &  0   & 0 \\
1  &  -1 & 0 \\
0  &  1  & -1 \\
\end{array}
\right ].  
\end{align}

\subsection{State and fault estimation using $\ell_1$-norm optimization}\label{estimator_sec}
For each $k\in\mathbb{N}$, we aim at estimating the state $x(k)$ and the fault $f(k)$, based on the current output measurement $y(k)$. We start by formulating a \textit{centralized} solution, which estimates them by employing all output measurements as in \req{overall_output2} or \req{overall_output}. A distributed strategy will be formulated later in this paper. Let $\hat{x} (k)$, $\hat{f} (k)$, $k\in\mathbb{N}$ be the state and the fault signals that are estimated at $k$, respectively. 
In order to obtain these estimates, we solve the following optimization problem: 
\begin{align}
\hat{x}(k) &=\underset{x}{{\rm arg\ min}}\ \| x - A \hat{x} (k-1) - B u(k-1) \| _1 \notag \\
 &{\rm s.t.}\ \ y (k) = C x, \label{optimization_problem}
\end{align}
where $C = C_0$ if $a_1 (k) = 0$ and $C=C_1$ if $a_1 (k) = 1$. Without loss of generality, we set $\hat{x} (k-1) = 0$ for $k=0$. The estimation for $\hat{f} (k)$ is then computed as 
\begin{equation}\label{anomaly_estimation}
\hat{f} (k) = \hat{x} (k) - A \hat{x} (k-1) - Bu(k-1). 
\end{equation}
In general, the term $A \hat{x} (k-1) + Bu(k-1)$ is called \textit{apriori estimate} for $k$. Thus, the optimization problem in \req{optimization_problem} aims at finding the state that is the closest to the apriori estimate, subject to the constraint on the output measurement $y (k)$. Letting $z = x - A \hat{x} (k-1) - Bu(k-1)$, the problem in \req{optimization_problem} leads to $\underset{x}{{\min}}\ \|z\| _1$, subject to $\widetilde{y}(k) = C z$, 
where $\widetilde{y}(k) = y (k)- C (A\hat{x}(k-1) + B u(k-1))$, which is indeed a well-known Basis Pursuit (BP) \cite{basispursuit2}. Hence, various numerical solvers can be applied to obtain the solution, both in a centralized manner \cite{basispursuit3}, and distributed manner \cite{basispursuit} (see also \rsec{distributed_section} for the distributed formulation). 

The optimization problem in \req{optimization_problem} is indeed relevant to the existing estimation methodologies. 
To see this, consider the following \textit{weighted $\ell_2$-norm} optimization problem: 
\begin{equation}\label{l2optimization}
\underset{x}{\min}\ \| x - A \hat{x} (k-1) - Bu(k-1) \|^2 _{P^{-1}} + \| y (k) - C x \|^2 _{V^{-1}}, 
\end{equation}
where $V, P$ are the positive definite matrices. 
Then, the solution to the above optimization problem is given by 
\begin{align}
\hat{x} (k) = &A \hat{x}^- (k-1) \notag \\ 
 &\!\!\!\!\!\!\!\!\! + A P C^\mathsf{T} (V + C P C^\mathsf{T})^{-1} (y(k) - C A \hat{x}^- (k-1)), \label{kalman_update} 
\end{align}
where $\hat{x}^{-} (k) = A \hat{x} (k-1) + Bu(k-1)$, which is indeed the (steady state) Kalman filter update. Thus, the optimization problem in \req{optimization_problem} is, roughly speaking, a modified version of the Kalman filter by utilizing the $\ell_1$-norm in the objective function. {The utilization of the $\ell_1$-norm is motivated by the fact that it potentially provides a better estimation accuracy than the Kalman filter ($\ell_2$-norm optimization), if the number of faulty nodes is small enough; namely, if $f(k)$ has a \textit{sparse} structure that contains many zero components. This benefit is indeed illustrated by the following example. } 

\smallskip
\smallskip
\noindent
\textit{ (Example) } : 
Suppose that ${\cal V} = \{1, 2, 3\}$ and $\{1, 2\}, \{2, 3 \} \in {\cal E}$, and the dynamics of node $i$ is scalar and given by $x_i (k+1) = x_i (k) + f_i (k)$, $\forall i\in\{1, 2, 3\}$. As with the example in \rsec{dynamics_sec}, we have $A = I_3$, $B=0$, and $C_0$, $C_1$ are given by \req{c_example}, and the measurement output is given by \req{overall_output2} and \req{overall_output}. 
Assume that $a_1 (k) = 1$, $\forall k \in [0, 19]$ and $a_1 (k) = 0$, $\forall k \in [20, 40]$. 
The initial states are $x_0 = [2; 4; 6] \in \mathbb{R}^3$ and the evolution of the states are illustrated in \rfig{example} (blue solid lines).  As shown in the figure, it is assumed that the fault occurs at $k=30$ for $x_1 (k)$ with $f_1 (k) = -3$. 

The estimation results by applying \req{optimization_problem} and the Kalman filter ($\ell_2$-norm optimization) are also illustrated in \rfig{example} (left figure). During the Kalman filter implementation, we simply use \req{l2optimization} with the parameters given by $P = 0.0001 I_3$ and $V = 0.0001 I_3$ (if $C=C_1$), $V = 0.0001 I_2$ (if $C=C_0$). While both approaches provide somewhat good estimates for a while, it is shown that the $\ell_1$-norm optimization approach dramatically overcomes the Kalman filter when the fault occurs at $k=30$. 
 Intuitively, this is due to the fact that the proposed approach utilizes the $\ell_1$-norm objective function; if $\hat{x}(k-1) \approx x (k-1)$, which may hold in this example since the system is observable for a while,  we then have $x (k) - A \hat{x} (k-1) - B u(k-1) \approx f(k)$, and solving \req{optimization_problem} may lead to be a good estimator as $f(k)$ is \textit{sparse} with only the first component taking non-zero values. 

Although the $\ell_1$-norm optimization approach becomes a good estimator in the above example, it leads to a failure when more nodes are faulty. To illustrate this, the right figure of \rfig{example} shows the result when both $x _1 $ and $x_2$ are faulty concurrently at $k=30$. Clearly, the $\ell_1$-optimization approach fails to reconstruct states and faults, which may imply that the estimation performance depends on the number of faulty nodes. Additionally, since $\hat{x}(k-1)$ is utilized for the estimator in \req{optimization_problem}, the estimation accuracy may also depend on the estimation error at the previous time; if these errors are too large, the estimation accuracy may be worse. Motivated by the above observations, in this paper we aim to provide answers to the following questions. 
\begin{itemize}
\item \textit{How many faulty nodes can be tolerated to provide correct estimation? 
\item How does the estimation error at $k-1$ affect the estimation error at $k$?} 
\end{itemize}

\begin{figure}[t]
  \begin{center}
   \includegraphics[width=8.5cm]{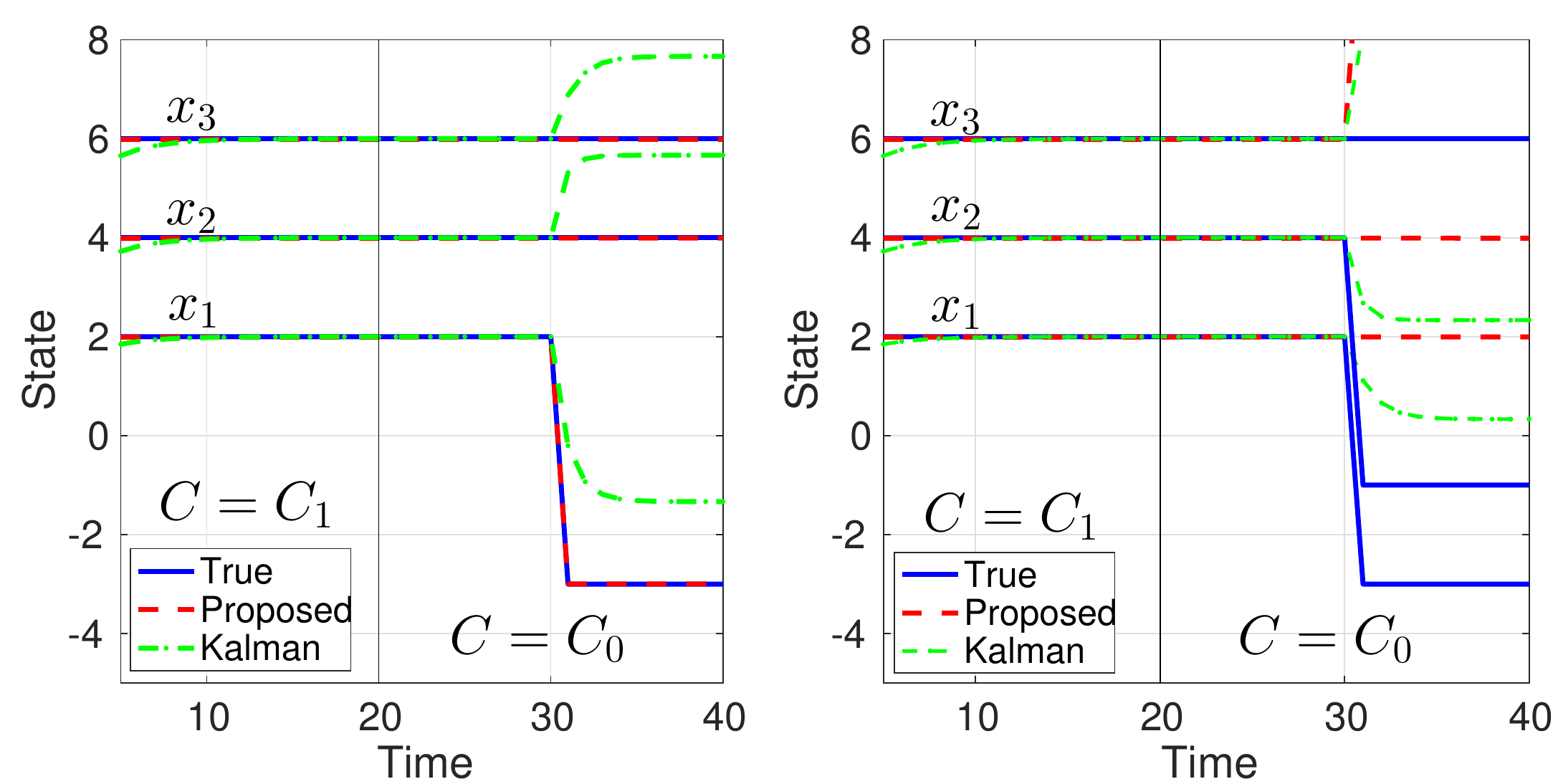}
   \caption{Example of applying \req{optimization_problem} and the Kalman filter when one node is faulty (left) and two nodes are faulty (right). The black vertical line indicates the time step when the system becomes un-observable ($k=20$).} 
   \label{example}
  \end{center}
  \vspace{-0.5cm}
 \end{figure}
 
\begin{myrem}[{On the extension to the case when measurement noise is present}]\label{noise_rem}
\normalfont
 {Suppose that measurement noise is present and we have $y(k) = C x(k) + w(k)$ ($C = C_0$ if the leader is non-active, and $C = C_1$ if the leader is active), where $w(k)$ represents measurement noise and is assumed to be norm-bounded as $\|w(k)\| \leq w_{\max}$, $\forall k\in\mathbb{N}$ for a given $w_{\max}\in\mathbb{R}_+$. In this case, it is possible to modify the estimator as follows: 
\begin{align}
\hat{x}(k) &=\underset{x}{{\rm arg\ min}}\ \| x - A \hat{x} (k-1) - B u(k-1) \| _1 \notag \\
 &{\rm s.t.}\ \ \| y (k) - C x \| \leq w_{\max}.  \label{optimization_problem_noise}
\end{align}

In this paper, we focus on providing a detailed theoretical analysis for the estimator of \req{optimization_problem}, rather than of \req{optimization_problem_noise}. 
However, in order to investigate the applicability of \req{optimization_problem_noise}, we also illustrate the corresponding estimation performance through a numerical simulation. For details, see \rsec{simulation_example}. \qedwhite } 
\end{myrem}

\section{Analysis of the estimator}
\label{analysis_no_error}
In this section we provide a quantitative analysis of the estimator as the answer to the first question considered in the previous section. 
Recall the optimization problem \req{optimization_problem}: 
\begin{equation}\label{optimization_problem_simple}
\underset{x}{{\min}}\ \| x - A \hat{x} (k-1) - Bu(k-1)\| _1 \ \ \ {\rm s.t.}\ \ y (k) = C x, 
\end{equation}
for $k\in\mathbb{N}$, where $C = C_0$ if $a_1 (k) = 0$ and $C=C_1$ if $a_1 (k) = 1$. 
Suppose that \ras{initial_assumption} holds ($a_1(0) = 1$, ${\cal I}_0 = \emptyset$), and consider the unknown sets of faulty nodes ${\cal I}_k \subseteq \{1, \ldots, M\} \cup \emptyset$, for $k\in\mathbb{N}\backslash \{0\}$, which means that $f_i (k) \neq 0$, $\forall i\in {\cal I}_k$ and $f_i (k) = 0$ otherwise. 
Let ${\cal F}_{{\cal I}_k} \subset \mathbb{R}^{nM}$ be given by 
\begin{align}
{\cal F}_{{\cal I}_k} = \{ f(k) \in \mathbb{R}^{nM}: f_i(k) \neq & 0, \forall i\in{\cal I}_k\ {\rm and}\ \notag \ \\ 
 &  \ f_i (k) = 0, \forall i\notin{\cal I}_k \}. 
\end{align}
That is, ${\cal F}_{{\cal I}_k}$ is the domain of all $f(k)$ when the set of faulty nodes is given by ${\cal I}_k$. Then, we can say that the estimator is \textit{correct} if the estimator in \req{optimization_problem_simple} yields the estimation that is identical to the actual one, i.e., for all $k\in \mathbb{N}\backslash\{0\}$, it follows that $\hat{x}(k) = x(k)$, $\hat{f}(k) = f(k)$, $\forall f(k) \in {\cal F}_{{\cal I}_k}$, $\forall u(k) \in \mathbb{R}^{mM}$. The following theorem presents a necessary and sufficient condition for this as the main result of this section:

\begin{mythm} \label{main_result_NSP}
\normalfont
Consider the estimator in \req{optimization_problem_simple}, and suppose that \ras{initial_assumption} holds. 
Then, the following statements are equivalent: 
\renewcommand{\labelenumi}{(\alph{enumi})}
\begin{enumerate}
\item $|{\cal I}_k| < M/2$, $\forall k\in\mathbb{N}\backslash \{0\}$. 
\item $\hat{x} (k) = x(k)$, $\hat{f} (k) = f(k)$, $\forall f(k) \in {\cal F}_{{\cal I}_k}$, $\forall u(k-1) \in \mathbb{R}^{mM}$, $\forall a_1(k) \in \{1, 0\}$, $\forall k\in\mathbb{N}\backslash \{0\}$. \qedwhite
\end{enumerate}
\end{mythm}
For the proof, see the \rapp{proof_appendix1}. 
In essence, \rthm{main_result_NSP} indicates that state and fault signals are correctly estimated for all $k\in\mathbb{N}\backslash \{0\}$ regardless of the values of $f(k)$, control inputs $u(k-1)$ and the leader's mode $a_1(k)$, if and only if the number of faulty nodes is less than half the total number of nodes for all $k\in\mathbb{N}\backslash \{0\}$. 
Indeed, this fundamental limitation is related to the ones presented in \cite{fawzi,shoukry2016,michelle}, where they showed that the number of attacks on sensors that can be tolerated to provide a correct estimation cannot exceed half the total number of sensors. 

\section{Analysis of estimator in the presence of estimation errors at previous time}\label{analysis_estimation_error}
In the previous section, we provide a necessary and sufficient condition such that the estimator provides correct estimations. 
In essence, correct estimations for $k \in \mathbb{N}\backslash\{0\}$ are ensured if and only if the estimation error at $k-1$ is absent, i.e., $x(k-1) = \hat{x}(k-1)$ (see \rlem{main_result_lem} in \rapp{proof_appendix1}). 
However, estimation errors at $k-1$ should arise due to, e.g., numerical errors of solving \req{optimization_problem}, imperfect knowledge of the initial state, the presence of system noise, and so on. 
In this section, therefore, we analyze how the estimation error at $k-1$ affects the estimation error at $k$ as the answer to the second question considered in \rsec{estimator_sec}. 
Moreover, we will also illustrate how the upper bound of the estimation error grows over time from the initial estimation error. 
The main result of this section is given as follows: 

\begin{mythm}\label{main_result_error}
\normalfont
For a given $k\in\mathbb{N}\backslash \{0\}$ and $d_{\max} \in \mathbb{R}_+$, suppose that 
$\|\hat{x} (k-1) - x(k-1)\|_1 \leq d_{\max}$ holds, and let ${\cal I}_k$ be any set of faulty nodes with $|{\cal I}_k| < M/2$. Then, for every $f(k) \in {\cal F}_{{\cal I}_k}$, $u(k) \in \mathbb{R}^{mM}$ and $a_1(k) \in \{0, 1\}$, it holds that 
\begin{equation}\label{error_bound}
\|f(k) - \hat{f} (k)\|_1 \leq \cfrac{2(M-|{\cal I}_k|)}{M-2|{\cal I}_k|}\ \eta  d _{\max}, 
\end{equation} 
where $\eta = \sum^{nM} _{i=1} \sum^{nM} _{j=1} |A^{(i, j)}|$. \qedwhite
\end{mythm}
For the proof, see the \rapp{proof_appendix2}. \rthm{main_result_error} states that, if the number of faulty nodes is smaller than half the total number of nodes, the estimation error of faults is bounded at most proportional to the estimation error of states occurred at the previous time. Moreover, since $\frac{2(M-|{\cal I}' _k|)} {M- 2|{\cal I}' _k|} < \frac{2(M-|{\cal I} _k|)} {M- 2|{\cal I}_k|}$ for any ${\cal I}' _k$, ${\cal I}_k$ with $|{\cal I}' _k| < |{\cal I} _k| < M/2$, the estimation error may be larger as the number of faulty nodes becomes larger. 

Based on \rthm{main_result_error}, we can obtain how the estimation error of states grows over time depending on whether the leader is active or non-active. To see this, suppose that 
$\|x(0) - \hat{x}(0)\|_1 \leq d(0)$, where ${d}(0) \in \mathbb{R}_+$ represents the upper bound of the estimation error for the initial time. Such estimation error may arise due to that: (i) the leader is in non-active mode at $k=0$ (i.e., \ras{initial_assumption} does not hold) and the initial state is unknown, but we know that the estimation error is bounded as $\|x(0) - \hat{x}(0)\|_1 \leq d(0)$; (ii) the leader is in active mode at $k=0$ (i.e., \ras{initial_assumption} holds), but some numerical errors arise in solving \req{optimization_problem}. 
Moreover, suppose that $|{\cal I}_1| < M/2$. 
If the leader is non-active at $k=1$, it then follows that 
\begin{align}
\|x(1) - \hat{x}(1)\|_1 & \leq \|A (x(0) - \hat{x}(0))\|_1 + \|f(1) - \hat{f}(1)\|_1 \notag \\
                            &\leq \cfrac{3M-4|{\cal I}_1|} {M- 2|{\cal I}_1|}\ \eta {d}(0), \label{next_error}
\end{align}
where we used $\|A(x(0) - \hat{x}(0))\|_1 \leq \eta {d}(0)$ and $\|f(1) - \hat{f}(1)\|_1 \leq \frac{2(M-|{\cal I}_1|)}{M-2|{\cal I}_1|} \eta {d}(0)$. 
On the other hand, if the leader is active at $k=1$, it follows that $C = C_1$ has a trivial kernel and the solution of $x$ satisfying $y(k) = Cx$ is uniquely determined (for details, see the proof of \rlem{main_result_lem} in \rapp{proof_appendix1}). 
Although this implies that $\hat{x}(1) = x(1)$ by solving \req{optimization_problem}, we consider here a more practical situation that $\|x(1) - \hat{x}(1)\|_1 \leq \bar{d}$, where $\bar{d} \in \mathbb{R}_+$ represents the upper bound of the numerical error in solving \req{optimization_problem}. 
Hence, by applying the above for $k=0, 1, 2, \ldots$, we obtain $\|x(k) - \hat{x}(k)\|_1 \leq d(k)$, where $d(k) \in \mathbb{R}_+$, $k \in \mathbb{N}$ are recursively given by
\begin{numcases}
{d(k)=} 
 \cfrac{3M-4|{\cal I}_k|} {M- 2|{\cal I}_k|} \eta d (k-1), \ \ {\rm if}\ a(k) = 0, \label{active_d} \\ 
 \bar{d} \ \ \ \ \ \ \ \ \ \ \ \ \ \ \ \ \ \ \  \ \ \ \ \ \ \ \ \ {\rm if}\ a(k) = 1, 
\end{numcases}
for all $k\in\mathbb{N}\backslash \{0\}$, where we assume that $|{\cal I}_k|<M/2$, $\forall k\in\mathbb{N}\backslash \{0\}$. The above relation implies that, the upper bound of the estimation error is \textit{reset} to $\bar{d}$ if the leader gets into the active mode, while on the other hand the error may diverge as long as the leader is in non-active mode. 

\begin{myrem}[On the effect of system noise]
\normalfont
 {The analysis can be extended to the case when (bounded) system noise is added to the system. 
To illustrate this, suppose that the dynamics is given by $x(k) = A x(k-1) + Bu(k-1) + v(k-1) + f(k)$, where $v$ is the system noise and is assumed to be bounded as $\|v(k)\|_1 \leq v_{\max}$, $\forall k\in\mathbb{N}$ for a given $v_{\max} > 0$. Then, \req{next_error} becomes $\|x(1) - \hat{x}(1)\|_1 \leq \|A (x(0) - \hat{x}(0))\|_1 + \|f(1) - \hat{f}(1)\|_1 + \|v(1)\|_1 \leq \frac{3M-4|{\cal I}_1|} {M- 2|{\cal I}_1|}\ \eta {d}(0) + v_{\max}$. 
Hence, the corresponding $d(k)$ is given by replacing the right hand side of \req{active_d} with $\frac{3M-4|{\cal I}_k|} {M- 2|{\cal I}_k|} \eta d (k-1) + v_{\max}$. \qedwhite }
\end{myrem}


\section{Distributed implementation}\label{distributed_section}
In this section we provide an algorithm that yields the solution to \req{optimization_problem} in a \textit{distributed} manner, which means that each node is able to estimate states and faults for all nodes by coordinating only with its neighbors. The approach follows the distributed basis pursuit (DBP) algorithm \cite{basispursuit}, which is a solver for $\ell_1$-norm optimization problem based on the alternating direction method of multipliers (ADMM). In this section, we provide an overview of the methodology presented in \cite{basispursuit} and apply to our fault and state estimation problem. 

Assume that the network graph ${\cal G} = ({\cal V}, {\cal E})$ is connected, bidirectional and bipartite: $\{i, j\} \in {\cal E}$ implies $\{j, i\} \in {\cal E}$ and the vertices can be decomposed into two independent sets ${\cal M}_1$, ${\cal M}_2$, i.e., every edge in ${\cal E}$ connects a vertex in ${\cal M}_1$ (or ${\cal M}_2$) to a vertex in ${\cal M}_2$ (or ${\cal M}_1$). 
In addition, assume that the control law for node $i$ is based on the relative state information from its neighbors, and all nodes share the same control law, i.e., 
\begin{equation}
u_i (k) = \rho (r_{{\cal N}_i} (k)),\ \forall i\in\{1, \ldots, M\}, 
\end{equation}
where $r_{{\cal N}_i} (k) = \{x_i(k) - x_j(k) : j \in {\cal N}_i \}$ and $\rho$ denotes the control law. 
The control law that collects all nodes is denoted as $u(k) = \kappa (x(k))$, where $\kappa = [\rho ^\mathsf{T}, \rho ^\mathsf{T}, \ldots, \rho ^\mathsf{T}]$. The following assumption is further required: 
\begin{myas}\label{uisknown}
\normalfont
The matrix $A$ is known to all nodes. \qedwhite
\end{myas}
\noindent
Recall that $y(k)$ and $C_0, C_1$ can be partitioned by rows as 
\begin{equation}
{\small C_q =  \left [
\begin{array}{c}
 C_{1q}  \\
 C_2 \\
 \vdots \\
 C_M
\end{array}
\right ],\ \ \ \ 
y(k) =  \left [
\begin{array}{c}
 y_1(k)  \\
 y_2(k) \\
 \vdots \\
 y_M(k)
\end{array}
\right ],}
\end{equation}
for $q \in \{0, 1\}$, where $y_i$, $i\in\{1, \ldots, M\}$ is the output measurements that are available to node $i$. To employ the approach presented in \cite{basispursuit} for the row partition case, let $\chi_i \in \mathbb{R}^{nM}$, $i\in \{1, \ldots, M\}$ be the decision variable that will be optimized by node $i$. Roughly speaking, $\chi_i (k)$ represents the estimated state of $x(k)$ obtained by node $i$. 
The optimization problem in \req{optimization_problem} is then re-formulated as follows: 
\begin{align}
  \min \ \cfrac{1}{M} & \sum^{M} _{i=1} \| \chi_i - A \hat{\chi}_i  (k-1) - B \kappa (\hat{\chi}_i  (k-1)) \| _1 \notag \\
{\rm subject\ to:\ }&\ y_1 (k) = C_{1q} \chi_1, \notag \\
&\ y_i (k) = C_i \chi_i,\ \  \forall i\in\{2, \ldots, M\} \notag \\  
&\ \chi_i = \chi_j, \ \ \forall \{i, j\} \in {\cal E}, \label{constraint_edge} 
\end{align}
where $q = 0$ (resp. $q=1$) if $a_1 (k) = 0$ (resp. $a_1 (k) = 1$), and $\hat{\chi}_i  (k-1) \in \mathbb{R}^{nM}$ is the estimated state obtained by node $i$ at $k-1$. Let $\hat{\chi}_i (k) \in \mathbb{R}^{nM}$, $i\in \{1, \ldots, M\}$ be the solution to the optimization problem in \req{constraint_edge}. 
From the connectivity of the graph and the constraint in \req{constraint_edge}, global consistency for the optimal solution is satisfied, i.e., $\hat{\chi}_1 (k) = \hat{\chi}_2 (k) = \cdots = \hat{\chi}_M (k)$. 
Using the independent sets ${\cal M}_1$, ${\cal M}_2$, the optimization problem in \req{constraint_edge} can be rewritten as 
\begin{align}
  \min\ \cfrac{1}{M} \sum_{i\in{\cal M}_1} &\| \chi_i - A \hat{\chi}_i  (k-1)- B \kappa (\hat{\chi}_i  (k-1)) \| _1 \notag \\ 
  &\! \! \! \! \! \! \! \! \!+ \cfrac{1}{M}\sum_{i\in{\cal M}_2} \| \chi_i - A \hat{\chi}_i  (k-1) - B \kappa (\hat{\chi}_i  (k-1)) \| _1 \notag \\
{\rm subject\ to:\ } &\ y_1 (k) = C_{1q} \chi_1, \notag \\
&\ y_i (k) = C_i \chi_i,\ \  \forall i\in\{2, \ldots, M\}, \notag \\  
&\ (E^\mathsf{T} _1 \otimes I_{nM}) \bar{\chi}_1 +  (E^\mathsf{T} _2 \otimes I_{nM})\bar{\chi}_2 = 0,  \label{constraint_edge2} 
\end{align}
where $\bar{\chi}_1 \in \mathbb{R}^{nM|{\cal M}_1|}$ and $\bar{\chi}_2\in \mathbb{R}^{nM|{\cal M}_2|}$ are the vectors collecting $\chi_i$ for all $i\in{\cal M}_1$ and $i\in{\cal M}_2$, respectively, $I_{nM}$ denotes the $nM \times nM$ identity matrix, $\otimes$ denotes the Kronecker product, and $E_1$, $E_2$ are the matrices such that $E = [E^\mathsf{T} _1, E^\mathsf{T} _2]^\mathsf{T}$ is the incidence matrix of the graph ${\cal G}$. 
The structure of the optimization problem in \req{constraint_edge2} is now suited to apply the ADMM algorithm in a distributed manner; by defining the augmented Lagrangian ${\cal L}(\bar{\chi}_1, \bar{\chi}_2, \lambda)$, it can be verified that minimizing ${\cal L}$ with respect to $\bar{\chi}_1$ and $\bar{\chi}_2$ can be executed parallelly among the nodes $i\in {\cal M}_1$ and $i\in {\cal M}_2$, respectively, and the dual variables can be also updated parallelly for all nodes, see \cite{basispursuit} for a detailed derivation. The overall distributed algorithm is shown in \ralg{overall_alg}. 
\begin{algorithm}[t]\label{overall_alg}
\small{ 
    Set $k=0$ and $\hat{\chi}_i (k-1) = 0$, $\forall i\in\{1, \ldots, M\}$; (initialization) \\
    \For{ each $k\in\mathbb{N}$ }{
    Set $\mu^1 _i = \chi^1 _1 = 0$, $\forall i\in \{1, \ldots, M\}$ and $L = 1$;  \\
    \While {$L < L_{\max}$} {
    		\For { all $i \in {\cal M}_1$ } {
    			Set $v^L _i = \mu^L _i - \zeta \sum_{j\in{\cal N}_i} \chi^L _j$ and solve the following problem: \label{step1}
    			\begin{align}
    			\hspace{-0.5cm}
    		 &\underset{\chi_i}{\min}\ \frac{1}{M} \| \widetilde{\chi}_i \|_1 + {(v^L _i)}^\mathsf{T} {\chi}_i + \frac{\nu_i \zeta}{2} \|{\chi}_i\|^2 \notag \\
    			&{\rm s.t.} \ {y}_i = C_i {\chi}_i,\ ({y}_i = C_{iq} {\chi}_i \ {\rm if}\ i=1), \label{update1}
    			\end{align}
    			where $\widetilde{\chi}_i = \chi_i - A \hat{\chi}_i (k-1)- B \kappa (\hat{\chi}_i  (k-1))$
    			and let $\chi^{L+1} _i$ be the optimal solution to \req{update1}. 
    			
    			Send $\chi^{L+1}_i$ to all $i \in {\cal N}_i$; \label{step2}
    		 }
    		 \For { all $i \in {\cal M}_2$ }{
    		  Repeat the same procedure in \rline{step1}\ --\ \rline{step2}; 
    		 }
    		 \For  { all $i \in {\cal V}$ }{
    		 Set 
    		 \begin{equation}
    		 \mu^{L+1}_i = \mu^L _i + \zeta \sum_{j\in{\cal N}_i} (\chi^{L+1} _i - \chi^{L+1} _j) 
    		 \end{equation}
    		 }
    		 Set $L:= L+1$; 
    		 }
    		 
    		 For all $i \in \{1, \ldots, M\}$, set $\hat{\chi}_i (k) = {\chi}^L _i$ and $\hat{\varsigma}_i (k) = \hat{\chi}_i (k) - A\hat{\chi}_i (k-1) - B \kappa (\hat{\chi}_i  (k-1))$; 
    }
    }
    \caption{Distributed implementation.} 
\end{algorithm}
In the algorithm, the variable ${\chi}_i$ is updated according to \req{update1}, where $\nu_i >0$ is the $i$-th diagonal element of $E_1E^\mathsf{T} _1$, and $\zeta>0$ is a constant parameter for the quadratic penalty added in the objective function to define the augmented Lagrangian. By \ras{uisknown}, the updates in \req{update1} can be executed in parallel for all nodes in ${\cal M}_1$. Once all updates in $i\in{\cal M}_1$ are done, these are transmitted to its neighbors. Similarly, the updates for all nodes in ${\cal M}_2$ are executed in parallel and are transmitted to their neighbors. Afterwards, the variables $\mu_i$ relating to the dual variables are updated in parallel for all nodes in $i\in{\cal V}$. The above procedure is iterated until the number of iterations reaches a prescribed threshold $L_{\max}$, and obtain the estimates of states and faults that are denoted as $\hat{\chi}_i (k)$ and $\hat{\varsigma}_i (k)$, respectively. 

Based on the above algorithm, the following convergence property is satisfied:
\begin{mylem}\label{convergence_result}
\normalfont
For a given $k\in\mathbb{N}\backslash \{0\}$ and $\bar{\chi} \in \mathbb{R}^{nM}$, let $\hat{x}(k)$ be the solution to \req{optimization_problem} with $\hat{x} (k-1) = \bar{\chi}$ and let $\hat{\chi}_i (k)$, $i\in\{1, \ldots, M\}$ be the estimated states by applying Algorithm~1 at $k$ with $\hat{\chi}_{i}(k-1) = \bar{\chi}$, $\forall i\in\{1, \ldots, M\}$. Then, it follows that $\hat{x} (k) = \hat{\chi}_i (k)$, $\forall i\in\{1, \ldots, M\}$ as $L_{\max} \rightarrow \infty$.  \qedwhite 


\end{mylem}

Roughly speaking, \rlem{convergence_result} shows that the state estimate at $k$ obtained by applying \ralg{overall_alg} converges to a solution to \req{optimization_problem}. 
The result follows exactly the same line as \cite{basispursuit} and is omitted in this paper. Recalling that $\hat{x}(k-1) = 0$ for $k=0$ (see \rsec{estimator_sec}), we have $\hat{x} (k-1) = \hat{\chi}_i (k-1)$, $\forall i\in\{1, \ldots, M\}$ for $k=0$. Thus, the following is immediate by applying \rlem{convergence_result} for all $k\in\mathbb{N}$:

\begin{mythm}
\normalfont
Suppose that \ralg{overall_alg} is implemented for all $k$. Then, it follows that $\hat{x} (k) = \hat{\chi}_i (k)$, $\forall i\in\{1, \ldots, M\}$, $\forall k\in\mathbb{N}$ as $L_{\max} \rightarrow \infty$. \qedwhite 
\end{mythm}

\begin{myrem}
\normalfont
 {By extending Algorithm~1, it is possible to incorporate attack (fault) detection and isolation strategies in several ways. For example, one may set a certain threshold $\varepsilon > 0$, and agent $i$ estimates that the neighboring agent $j$ ($j \in {\cal N}_i$) is faulty if $\| (\hat{\varsigma}_i)_{F_j} (k) \| > \varepsilon$, where $F_j = \{(j-1)n +1, (j-1)n+2, \ldots, jn \}$ (i.e., $(\hat{\varsigma}_i)_{F_j} (k)$ represents the estimation of $f_j (k)$ computed by agent $i$). Once the fault has been detected, agent $i$ removes the connection to agent $j$ in order to isolate it from the neighbors.} \qedwhite 
\end{myrem}

\section{Illustrative example}\label{simulation_example}
In this section we provide an illustrative example to
validate our estimation scheme. 
We consider a multi-vehicle system with $9$ vehicles (nodes), which moves in two-dimensional plane. Let $p_{i} = [p_{x_i}; p_{y_i}] \in \mathbb{R}^2$, $v_{i} = [v_{x_i}; v_{y_i}] \in \mathbb{R}^2$, $i\in\{1, \ldots, 9\}$ be the position and the velocity of vehicle $i$, respectively. The dynamics of vehicle $i$ is assumed to be a double-integrator: 
\begin{equation} \label{sim_system}
\dot{{x}}_i  =  \left [
\begin{array}{cccc}
0  &  0 & 1 & 0  \\
0 &  0  & 0 & 1  \\
0 & 0 & 0& 0 \\
0 & 0 & 0 & 0
\end{array}
\right ] x_i + \left [
\begin{array}{cc}
0 & 0 \\
0 & 0 \\
1 & 0 \\ 
0 & 1 \\
\end{array}
\right ] u_i + f_i, 
\end{equation}
where $x_i = [p_i ; v_i] \in \mathbb{R}^4$ and $u_i = [u_{x_i}; u_{y_i}] \in \mathbb{R}^2$ is the control input, and $f_i \in \mathbb{R}^4$ is the fault signal. 
We discretize the system in \req{sim_system} under a zero-order hold with $0.05s$ sampling time interval to obtain $x_i (k) = A_i x_i(k-1) + B_i u_i(k-1) + f_i(k)$, where $A_i$ and $B_i$ are the appropriate matrices. We assume that only velocity states are subject to faults, i.e., $f_i(k) = [0\ 0\ f_{x_i}\ f_{y_i}]^\mathsf{T}$. Since we have $\dot{v}_{x_i} = u_{x_i} + f_{x_i}$ and $\dot{v}_{y_i} = u_{y_i} + f_{y_i}$, one can also see that this situation is the case when actuators are subject to faults by physical faults or by malicious attackers. 

\begin{figure}[t]
  \begin{center}
   \includegraphics[width=2.5cm]{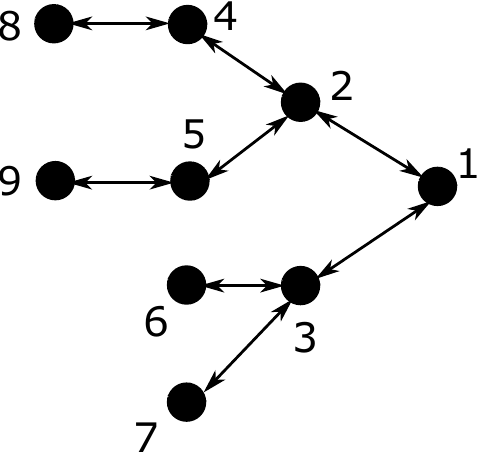}
   \caption{Graph structure of agents. Each node represents a vehicle and each edge represents a communication-and-sensing capacity between the corresponding vehicles. } 
   \label{sim_illust}
   \vspace{-0.5cm}
  \end{center}
 \end{figure}

The network graph ${\cal G} = ({\cal V}, {\cal E})$ indicating communication and sensing capacities among the vehicles is shown in \rfig{sim_illust}. 
We assume that there exist two environments; \textit{outdoor} and \textit{indoor environment}. In the outdoor environment, the leader is able to communicate with GPS satellites to obtain the state information $x_1$, which allows the leader to be in active mode. 
In the indoor environment, the leader is not able to utilize the GPS due to signal loss by presence of walls in the building, which is thus forced to be in non-active mode. The leader initially starts in the outdoor environment to ensure \ras{initial_assumption} and it enters the indoor environment afterwards. 
More specifically, we assume that active/non-active mode is given by
\begin{numcases}
{a_1 (k) =}
0, \ \ {\rm\ if}\ k\in [100, 300], \\
1, \ \ {\rm\ otherwise}. 
\end{numcases}
Note that for the overall system derived in \req{overall_state}, the pair $(A, C_1)$ is observable, while $(A, C_0)$ is \textit{not} observable. This means that the observability is lost while the leader is non-active during $k\in[100, 300]$. 

The control law for the leader is given by $u_{x_1} = - v_{x_1} + 1$, $u_{y_1} = 0$, so that it moves along the $x$-axis with a constant velocity $v_{x_1} = 1$. All the other vehicles aim to follow the leader by interacting with their neighbors. Specifically, the control law for the followers are given by $u_{x_i} = \sum_{j\in{\cal N}_i} -c_1 (p_{x_i} - p_{x_j} + d_{x_{ij}})$, 
$u_{y_i} = \sum_{j\in{\cal N}_i} -c_2 (p_{y_i} - p_{y_j} + d_{y_{ij}})$, 
where $c_1, c_2 >0$ are given constant gains and $[d_{x_{ij}}; d_{y_{ij}}]$ is the desired distance vector from $i$ to $j$.
\begin{figure}[t]
   \centering
    \subfigure[Faulty signals $f_{xi}$ and $f_{yi}$.]{
      {\includegraphics[width=7.5cm]{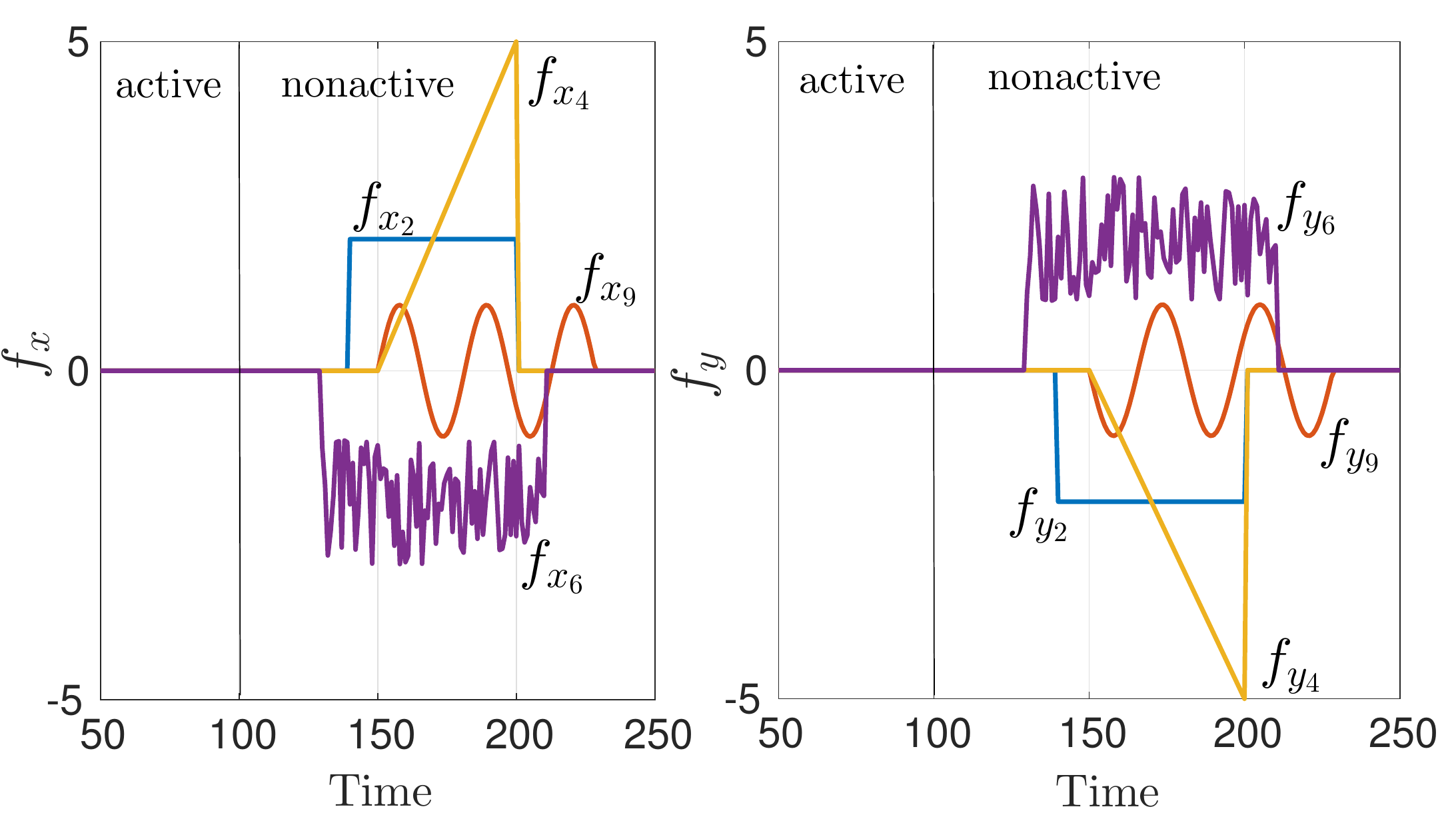}} }
    \subfigure[Trajectories of $p_{x_i}$.]{
      {\includegraphics[width=7.5cm]{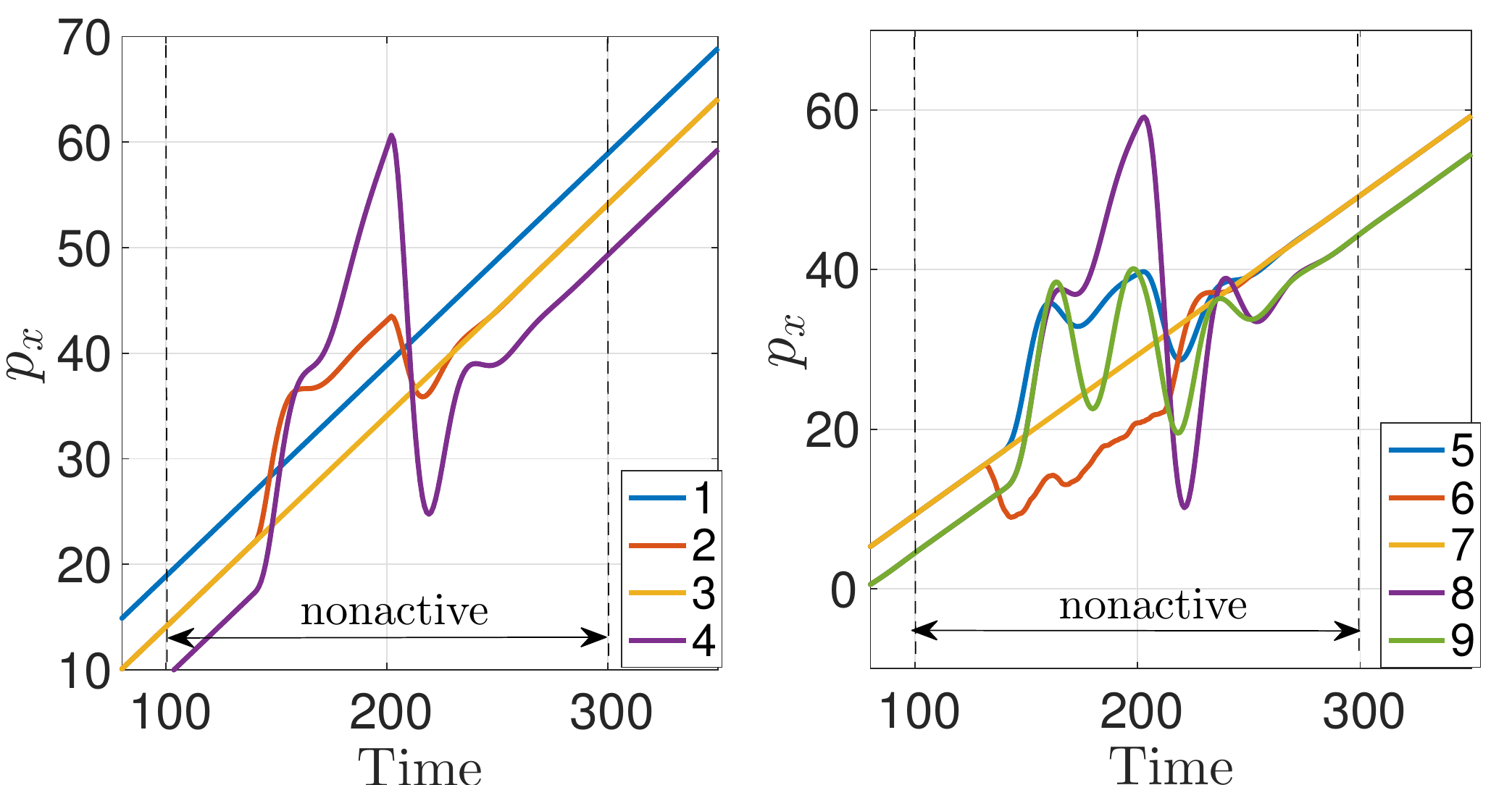}} }
    \caption{True trajectories of fault signals $f_{xi}$, $f_{yi}$ injected to the vehicles and the states of $p_{x_i}$.} \label{true_trajectory}
    \vspace{-0.5cm}
\end{figure}
\begin{figure}[t]
   \centering
    \subfigure[Estimated fault signals $f_{xi}$ and $f_{yi}$.]{
      {\includegraphics[width=7.5cm]{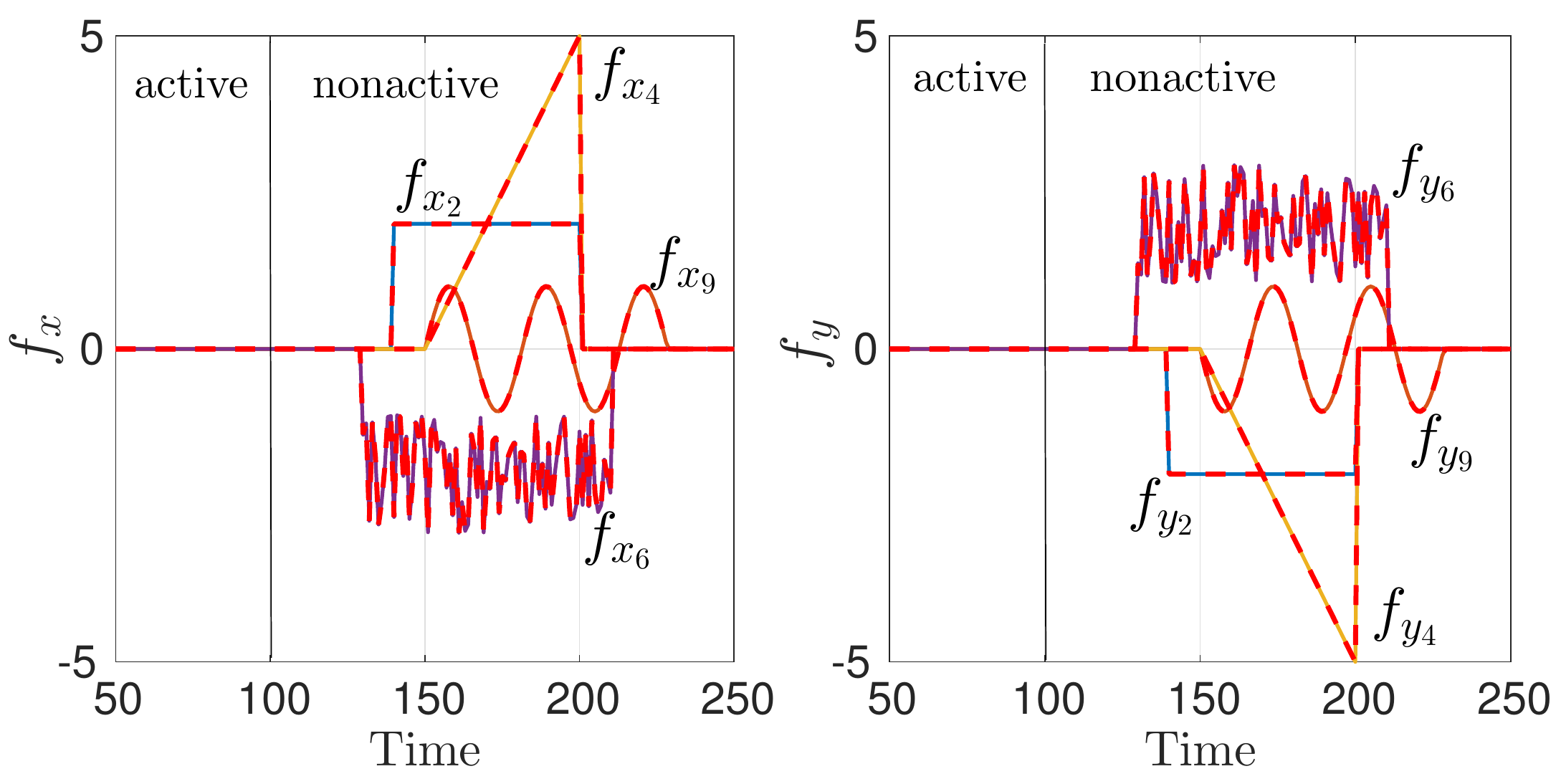}} \label{anomaly_signal_estimate}}
    \subfigure[Estimated trajectories of $p_{x_i}$.]{
      {\includegraphics[width=7.5cm]{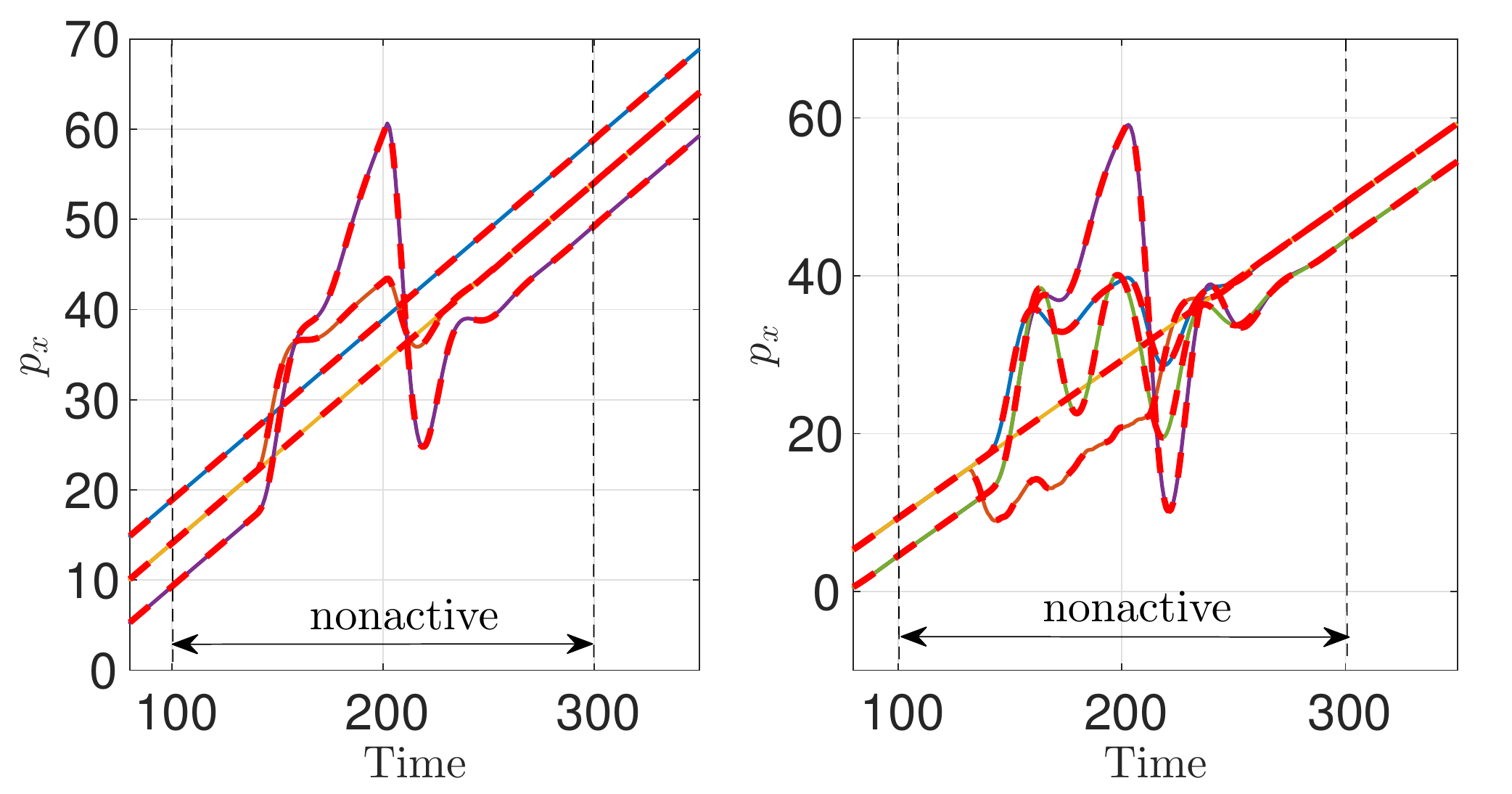}} \label{sample_traj_estimate}}
    \caption{Estimated fault signals and states of $\chi_1(k)$ by applying Algorithm~1 (red dotted lines).} 
    \label{state_trajectory_estimate}
    \vspace{-0.5cm}
\end{figure}
\begin{figure}[t]
  \begin{center}
   \includegraphics[width=7.5cm]{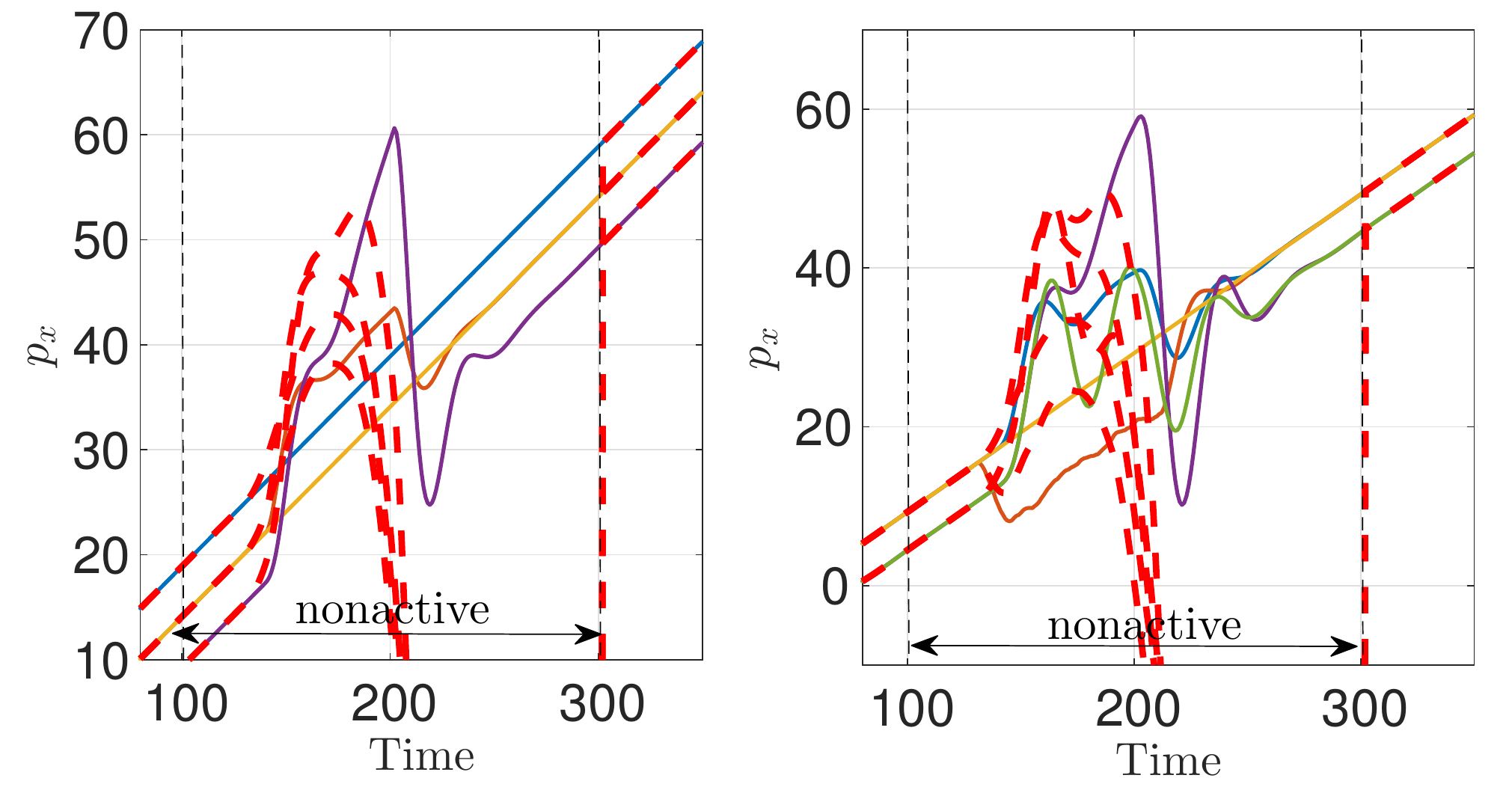}
   \caption{Estimated states by applying the Kalman filter (red dotted lines)} 
   \label{state_trajectory_estimate_kalman}
   \vspace{-0.5cm}
  \end{center}
 \end{figure}
\begin{figure}[t]
   \centering
    \subfigure[Estimation error of states against the time step.]{
      {\includegraphics[width=3.7cm]{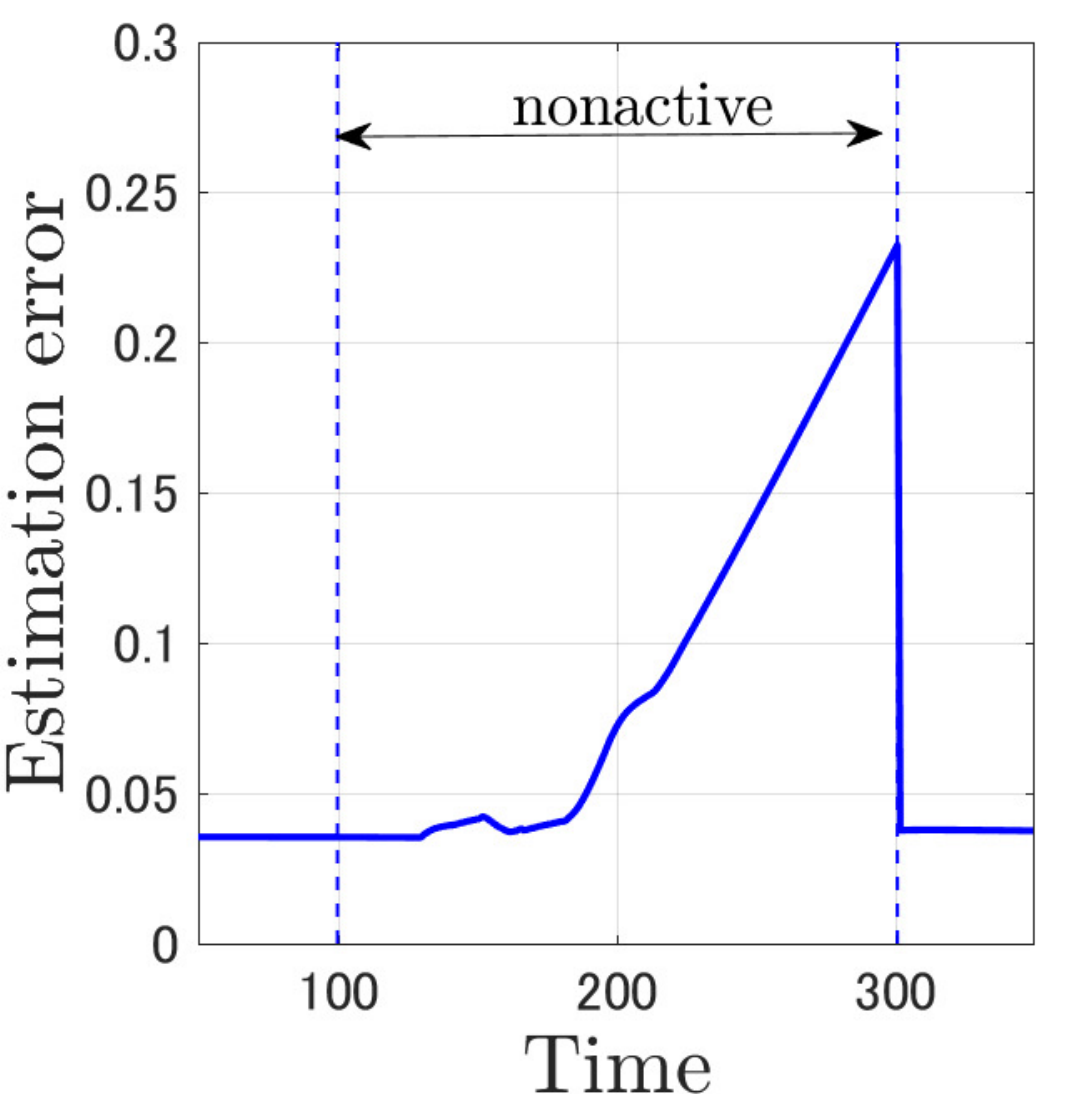}} \label{estimation_error_time}}
       \hspace{-0cm}
    \subfigure[Estimation error of states against the number of anomalous nodes.]{
      {\includegraphics[width=3.7cm]{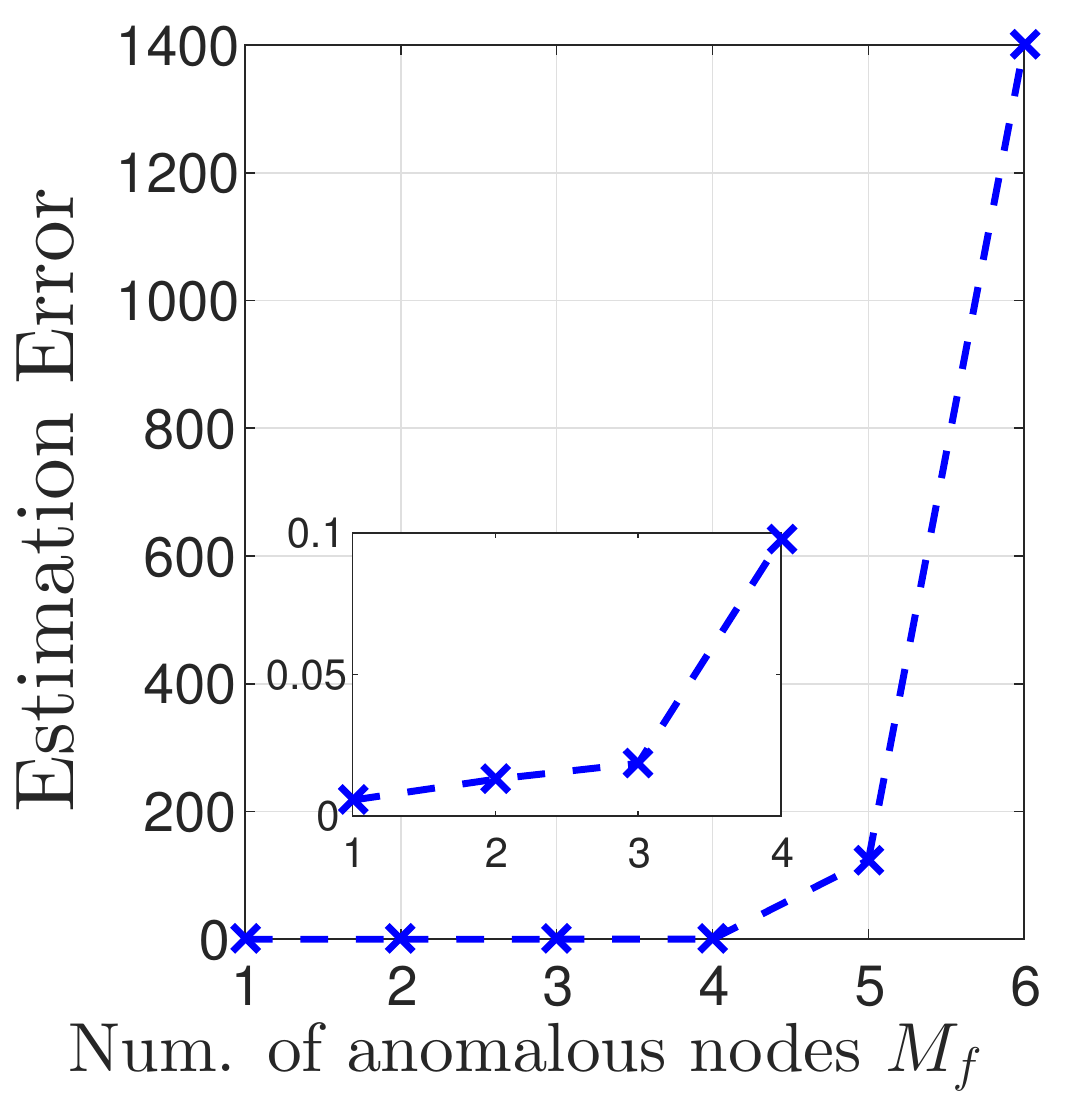}} \label{estimation_error}}
    \caption{Estimation error of states against the time step (left) and against the number of anomalous nodes (right). } 
    \vspace{-0.5cm}
\end{figure}
\rfig{true_trajectory} illustrates fault signals that are injected to the vehicles and the actual trajectories of $p_i$. As shown in the figure, it is assumed that the vehicle $2, 4, 6$ and $9$ are subject to faults with different shapes while the leader is the non-active mode.  As a consequence, the resulting trajectories provide somewhat faulty behaviors as shown in \rfig{true_trajectory}. 
 
\rfig{state_trajectory_estimate} illustrates the estimated states and faults by the leader (i.e., $\chi_1(k)$) by applying Algorithm~1. As shown in the figure, state and fault signals are appropriately estimated by applying the proposed approach. This is due to the fact that the number of faulty vehicles is given by $4 <9/2$, which satisfies the condition to guarantee correct estimation from \rthm{main_result_NSP}. 
 {For comparisons, we illustrate in \rfig{state_trajectory_estimate_kalman} the estimated trajectories by applying the Kalman filter ($\ell_2$-norm optimization). While implementing the Kalman filter, we utilize the update in \req{l2optimization} (or \req{kalman_update}) and the fault signals are estimated according to \req{anomaly_estimation}. From the figure, the estimation error {diverges} while the leader is non-active. The main reason for this is that the observability is lost when the leader is non-active, and the Kalman filter becomes vulnerable to the fault signals. Thus, the result shows that the proposed approach is more resilient against faulty signals than the Kalman filter.} 

While the estimated states in \rfig{state_trajectory_estimate} seem to exactly match the actual ones, the estimation error actually increases due to the numerical errors of solving the optimization problem. 
To see this, we illustrate in \rfig{estimation_error_time} the sequence of estimation error $\|x(k) - \hat{\chi_1}(k)\|_2$ for $k\in[50, 350]$. As shown in the figure, the error becomes larger as time evolves while the leader is in non-active mode, and it becomes much smaller as soon as it becomes active again at $k=301$. Indeed, this corresponds to the observation described in \rsec{analysis_estimation_error}; the upper bound of the error grows if the leader is non-active, while it is reset to an arbitrary small value when it becomes active. To see how the number of faulty nodes affect the estimation error, we illustrate in \rfig{estimation_error} the cumulative estimation error of states computed as $\sum^{300} _{k=101} \| \hat{\chi}_1(k) - x(k) \|_2$ with respect to the number of faulty vehicles $M_f  \in \{1, \ldots, 6\}$. For each $M_f \in \{1, \ldots, 6\}$, the set of faulty vehicles are selected as ${\cal I}_k = \{1, \ldots, M_f\}$ for all $k \in [100, 300]$, and $f_{x, i} (k)$, $i\in \{1, \ldots, M_f \}$ is assumed to be a random noise occurring with $f_{x, i} (k) \in [-10, 10]$. 
From the figure, it is shown that the estimation error dramatically increases when the number of anomalous vehicles is larger than $4$, which indeed violates the condition of $|{\cal I}_k| < M/2$ derived in \rthm{main_result_NSP}. Moreover, one can see from the enlarged view in the figure that the estimation error becomes larger even while $|{\cal I}_k| < M/2$. As is described in the observation in \rsec{analysis_estimation_error}, this is potentially because the upper bound of the estimation error  becomes larger as $|{\cal I}_k|$ becomes larger.

\begin{figure}[t]
  \begin{center}
   \includegraphics[width=8cm]{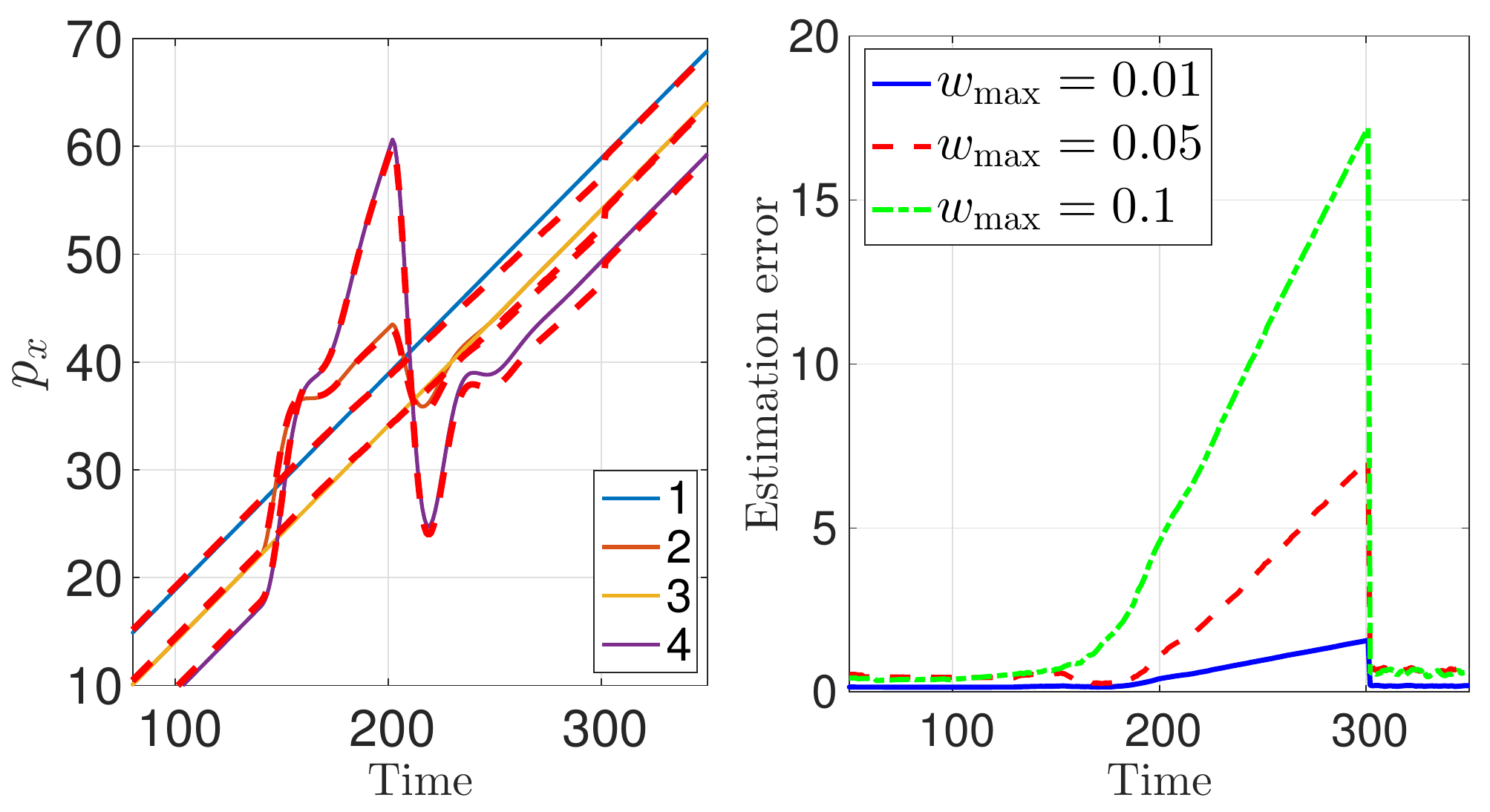}
   \caption{Estimated states by solving \req{optimization_problem_noise} with $w_{\max} = 0.05$ (left) and the sequence of the estimation error (right).} 
   \label{noise_error}
   \vspace{-0.6cm}
  \end{center}
 \end{figure}

 {Recall in \rrem{noise_rem} that we can potentially take measurement noise into account by solving \req{optimization_problem_noise}. 
To analyze the performance when measurement noise is present, we illustrate in \rfig{noise_error} the estimated states by solving \req{optimization_problem_noise} with $w_{\max} = 0.05$ (left) and the sequence of the estimation error $\|x(k) - \hat{x}(k)\|_2$ for $k\in[50, 350]$ with $w_{\max} = 0.01, 0.05, 0.1$ (right). 
The figure shows some robustness against measurement noise; although we did not provide any theoretical analysis, the estimation error seems to remain small if measurement noise is small enough. 
On the other hand, the figure demonstrates some practical limitations of our approach, since the estimation error grows over time when the leader is non-active, and the error becomes larger as $w_{\max}$ becomes larger. Hence, it is worth analyzing how the estimation error grows over time when measurement noise is present, and it is indeed one of the topics to be investigated in our future research.}

\section{Conclusions}
In this paper, we investigated a distributed state-and-fault estimation algorithm for multi-agent systems. 
The estimator employs an $\ell_1$-norm optimization problem, which is based on the signal recovery in compressive sampling. We then provided two quantitative analysis that characterize correctness of the proposed estimator. The results provide a fundamental limitation on the number of faulty nodes that can be tolerated to provide correct estimation, as well as how the estimation error is bounded and affected by the previous estimation error. Finally, a numerical example of multi-vehicle systems illustrated the effectiveness of the proposed approach.


\appendices
\section{Proof of \rthm{reconstruct}}
To prove sufficiency, suppose that $A \in \mathbb{R}^{m\times n}$ satisfies $T$-NSP for a given $T \subseteq \{1,\ldots, n\}$ and let $x_0 \in \mathbb{R}^n$ be an arbitrary $T$-sparse vector. Then, for any $v\in {\rm ker}(A)\backslash \{0\}$, it follows that 
\begin{align}
\|x_0 + v\|_1 - \|x_0\|_1 &= \sum^n _{i=1} |x^{(i)} _0 + v^{(i)} |  - \sum^n _{i=1} |x^{(i)} _0 | \notag \\
                             & = \sum_{i\in T} |x^{(i)} _0 + v^{(i)}|  + \sum_{i\in T^c} |v^{(i)}| - \sum _{i\in T} |x^{(i)} _0 | \notag  \\ 
                             & = \sum_{i\in T} (|x^{(i)} _0 + v^{(i)}|  - |x^{(i)} _0|) + \sum_{i\in T^c} |v^{(i)}| \notag \\
                             &\geq \sum_{i\in T^c} |v^{(i)}| - \sum_{i\in T} |v^{(i)}| \notag \\ 
                             &>0, \notag 
\end{align}
where $T^c = \{1, \ldots, n\}\backslash T$. For the second equality we used $x^{(i)} _0 = 0$, $\forall i\in T^c$ since $x_0$ is $T$-sparse, and for the third inequality we used the triangle inequality. For the last inequality, we used the $T$-NSP of $A$. This implies that among all $x$ satisfying $Ax_0 = Ax$, $\|x\|_1$ has a unique minimum when $x = x_0$. Hence, if $A$ satisfies $T$-NSP, every $T$-sparse vector $x_0$ is the unique solution to \req{l1_original}. 

To prove the necessity, suppose that every $T$-sparse vector $x_0$ is a unique solution to \req{l1_original}. 
For the sake of contradiction, assume that $A$ does \textit{not} satisfy $T$-NSP, i.e., there exists $v \in {\rm ker}(A)\backslash \{0\}$ such that $\|v_T\|_1 \geq \|v_{T^c}\|_1$ holds. Let $x_0$ be given by $x^{(i)} _0 = -v^{(i)}$, $\forall i \in T$ and $x^{(i)} _0 = 0$, $\forall i \in T^c$, which fulfills the $T$-sparsity assumption on $x_0$. Then, we have 
\begin{align}
\|x_0 + v \|_1 &= \sum _{i \in T} |x^{(i)} _0 + v^{(i)}|  + \sum _{i\in T^c } |x^{(i)} _0  + v^{(i)}| \notag \\ 
                 &= \sum _{i \in T} |-v^{(i)} + v^{(i)} | + \sum _{i\in T^c } |v^{(i)}| \notag \\ 
                 & = \sum _{i\in T^c } |v^{(i)}| \leq \sum _{i\in T } |v^{(i)}| \notag \\ 
                 & = \sum _{i\in T } |x^{(i)} _0| = \sum^n _{i=1 } |x^{(i)} _0| \notag \\
                 & = \|x_0 \|_1. \notag 
\end{align}
Thus, it follows that $\|x_0 + v \|_1 \leq \|x_0\|_1$. 
This implies that among all $x$ satisfying $A x_0 = Ax$, $\|x\|_1$ does not provide a unique minimum when $x = x_0$, which contradicts the fact that $x_0$ is a unique solution to \req{l1_original}. Therefore, if every $T$-sparse vector $x_0$ is a unique solution to \req{l1_original}, the matrix $A$ satisfies $T$-NSP. \qedwhite

\section{Proof of \rthm{main_result_NSP}}\label{proof_appendix1}
In order to provide the proof of \rthm{main_result_NSP}, we resort the following lemma:
\begin{mylem}\label{main_result_lem}
\normalfont
For a given $k\in\mathbb{N}\backslash \{0\}$ suppose that $\hat{x} (k-1) = x(k-1)$ holds. Then, the following statements are equivalent: 
\renewcommand{\labelenumi}{(\alph{enumi})}
\begin{enumerate}
\item $|{\cal I}_k| < M/2$. 
\item $\hat{x} (k) = x(k)$, $\hat{f} (k) = f(k)$, $\forall f(k) \in {\cal F}_{{\cal I}_k}$, $\forall u(k-1) \in \mathbb{R}^{mM}$, $\forall a_1(k) \in \{1, 0\}$. \qedwhite
\end{enumerate}
\end{mylem}

\rlem{main_result_lem} indicates that if there is no estimation error at $k-1$, a correct estimation is given at $k$ if and only if $|{\cal I}_k| < M/2$. 
Let us provide the proof of \rlem{main_result_lem}. 
To this end, it is required to take several steps to modify the optimization problem in \req{optimization_problem_simple}. 
Since $x (k-1) - \hat{x} (k-1) = 0$, we have $x(k) = A \hat{x} (k-1) + B u(k-1) + f(k)$ and so $y(k) = C(A \hat{x} (k-1) + B u(k-1) + f(k))$. Thus, the optimization problem becomes 
\begin{align}
&\underset{x}{{\min}}\ \| x - A \hat{x} (k-1) - Bu(k-1) \| _1 \\
 & {\rm s.t.,}\ \ C(A \hat{x} (k-1) + f(k) - x) = 0. 
\end{align}
Letting $z = x - A \hat{x} (k-1) - Bu(k-1)$ be the new decision variable to obtain
\begin{equation}\label{optimization_problem_z}
\underset{z}{{\min}}\ \| z \| _1 \ \ \ {\rm s.t.,}\ \ Cf(k) =  C z, 
\end{equation}
where $C = C_0$ if $a_1 (k) = 0$ and $C=C_1$ if $a_1 (k) = 1$. Let $\hat{z} (k)$ be the solution to \req{optimization_problem_z}. 
Since $z = x - A \hat{x} (k-1) - Bu(k-1)$ and $\hat{f}(k)$ is given by \req{anomaly_estimation}, we have $\hat{z} (k) = \hat{f} (k)$. 
Let $z = [z_1^\mathsf{T}, z_2^\mathsf{T}, \ldots, z_M^\mathsf{T}]^\mathsf{T}$, where $z_i \in\mathbb{R}^n$, $i\in \{1, \ldots, M\}$ represents the variable for node $i$. In the following, we rearrange the vectors $z$, $f(k)$ and the matrix $C$ as follows. Let $T_i \subset \{1, \ldots, nM\}$, $i\in \{1, 2, \ldots n\}$ be given by $T_i = \{i, i+n, i+ 2n, \ldots, i+ (M-1) n \}$, and rearrange $z$ as $z_p = [z_{p_1}^\mathsf{T}, z_{p_2}^\mathsf{T}, \ldots, z_{p_n}^\mathsf{T}]^\mathsf{T}$, where $z_{p_i} = z_{T_i}$, $i\in\{1, 2, \ldots, n\}$. 
For example, if $M=3$, $n = 2$ and
\begin{equation}
z = [\underbrace{z_{11}; z_{12}}_{z_1};\ \underbrace{z_{21}; z_{22}}_{z_2};\ \underbrace{z_{31}; z_{32}}_{z_3}] \in \mathbb{R}^6, 
\end{equation}
we have 
\begin{equation}
z_p = [\underbrace{z_{11}; z_{21}; z_{31}}_{z_{p_1}};\ \underbrace{z_{12}; z_{22}; z_{32}}_{z_{p_2}}]. 
\end{equation}
Roughly speaking, $z_{p_i}$ represents a vector collecting the $i$-th component of all the nodes. 
Similarly, rearrange $f(k)$ as $f_{p} = [f_{p_1}^\mathsf{T}, f_{p_2}^\mathsf{T}, \ldots, f_{p_n}^\mathsf{T}]^\mathsf{T}$, where $f_{p_i} = f_{T_i}$, $i\in\{1, 2, \ldots, n\}$. 
Note that due to the rearrangement of $f(k)$, $f_i (k) \neq 0$ implies that $(f_p)_{\tilde{T}_i} \neq 0$, where $\widetilde{T}_i = \{i, i+M, i+ 2M, \ldots, i+ (n-1)M\}$. 
For example, if 
\begin{equation}
f = [\underbrace{f_{11}; f_{12}}_{f_1};\ \underbrace{f_{21}; f_{22}}_{f_2};\ \underbrace{f_{31}; f_{32}}_{f_3}] \in \mathbb{R}^6, 
\end{equation}
we have 
\begin{equation}
f_p = [\underbrace{f_{11}; f_{21}; f_{31}}_{f_{p_1}};\ \underbrace{f_{12}; f_{22}; f_{32}}_{f_{p_2}}], 
\end{equation}
and $f_1 \neq 0$ implies $(f_p)_{\widetilde{T}_1} \neq 0$, where $\widetilde{T}_1 = \{i, i+M\} = \{1, 4\}$. 
Thus, letting $\widetilde{\cal F}_{{\cal I}_k}$ be given by 
\begin{align}
\widetilde{\cal F}_{{\cal I}_k} =  \{ f_p (k) \in \mathbb{R}^{nM}: (f_p & (k))_{\widetilde{T}_i} \neq 0, \forall i\in{\cal I}_k\ {\rm and}\ \notag \\
 & (f_p (k))_{\widetilde{T}_i} = 0, \forall i\notin{\cal I}_k \}, 
\end{align}
it follows that 
\begin{equation}
f(k) \in {\cal F}_{{\cal I}_k}\ \Longleftrightarrow\ f_p (k) \in \widetilde{\cal F}_{{\cal I}_k}. 
\end{equation}
Note that $f_p (k) \in \widetilde{\cal F}_{{\cal I}_k}$ implies that $f_p (k)$ is $\widetilde{T}$-sparse, where $\widetilde{T}= \cup_{i\in{\cal I}_k}\widetilde{T}_i$. 
Moreover, by rearranging the columns of $C$ according to the rearrangement for $z$ as well as suitably rearranging the rows, one can always construct the matrix $C_{p0}$ for the active mode ($a_1(k) = 1$) and $C_{p1}$ for the non-active mode ($a_1(k) = 0$) as 
\begin{align}
 C_{p0} &=  \left [
\begin{array}{cccc}
 D_1 &  0 & \cdots & 0 \\
 0 &  D_1 & \cdots & 0 \\
\vdots & \vdots  & \ddots & \vdots \\
 0   & 0 &   \cdots    & D_1
\end{array}
\right ], \label{cp} \\
C_{p1} &=  \left [
\begin{array}{cccc}
 D_2 &  0 & \cdots & 0 \\
 0 &  D_2 & \cdots & 0 \\
\vdots & \vdots  & \ddots & \vdots \\
 0   & 0 &   \cdots    & D_2
\end{array}
\right ], \label{cp2}
\end{align} 
where $D_1 = D^\mathsf{T}$ and $D_2 = [\mathsf{e} , D]^\mathsf{T}$ with $D$ being the incidence matrix of the graph ${\cal G}$ and $\mathsf{e} = [1, 0, 0, \cdots ,0]^\mathsf{T} \in \mathbb{R}^M$. 
For example, suppose that $C_0$, $C_1$ are given by
{ 
\begin{align}
{\small C_0 =  \left [
\begin{array}{cccc}
 1& 0 & -1 & 0 \\
 0& 1 & 0 & -1 
\end{array}
\right ], \ \ 
C_1 = \left [
\begin{array}{cccc}
1& 0 & 0 & 0 \\
0& 1 & 0 & 0 \\
 1& 0 & -1 & 0 \\
 0& 1 & 0 & -1 \\
\end{array}
\right ]}, \notag 
\end{align} 
which implies that $\{1, 2\} \in {\cal E}$ and the incidence matrix is given by $D = [1,  -1]^\mathsf{T}$. By rearranging the columns of the above matrices according to the one for $z$, we obtain 
\begin{align}
{\small C'_0 =  \left [
\begin{array}{cccc}
 1& -1 & 0 & 0 \\
 0& 0 & 1 & -1 
\end{array}
\right ], \ \ 
C'_1 = \left [
\begin{array}{cccc}
1& 0 & 0 & 0 \\
0& 0 & 1 & 0 \\
 1& -1 & 0 & 0 \\
 0& 0 & 1 & -1 \\
\end{array}
\right ],} \notag 
\end{align} 
By exchanging the second and the third rows of $C' _1$, we obtain 
\begin{align}
{\small C_{p0} =  \left [
\begin{array}{cccc}
 1& -1 & 0 & 0 \\
 0& 0 & 1 & -1 
\end{array}
\right ], \ \ 
C_{p1} = \left [
\begin{array}{cccc}
1& 0 & 0 & 0 \\
1& -1 & 0 & 0 \\
 0& 0 & 1 & 0 \\
 0& 0& 1 & -1 \\
\end{array}
\right ],} \notag 
\end{align} 
which can be checked to be of the form \req{cp} and \req{cp2}, respectively. Based on the above rearrangement, the optimization problem becomes 
\begin{equation}\label{optimization_problem_modified}
\underset{z_p}{{\min}}\ \| z_p \| _1 \ \ \ {\rm s.t.,}\ \ C_p f_p (k) = C_p z_p, 
\end{equation}
where $C_p = C_{p0}$ if $a_1(k) = 0$ and $C_p = C_{p1}$ if $a_1(k) = 1$. Let $\hat{f}_p (k)$ be the optimal solution of $z_p$ to \req{optimization_problem_modified}.
Note that the problem is invariant under the above rearrangements. That is, letting $\hat{z}(k) = \hat{f} (k)$ and $\hat{f}_p (k) =  [\hat{f}_{p_1}(k) ^\mathsf{T}$,$ \hat{f}_{p_2}(k)^\mathsf{T}$, $\ldots, \hat{f}_{p_n}(k)^\mathsf{T}]^\mathsf{T}$ be the solutions to \req{optimization_problem_z} and \req{optimization_problem_modified}, respectively, it follows that $\hat{f}_{p_i}(k) = \hat{f}_{T_i} (k)$, $\forall i\in\{1, \ldots, M\}$, i.e., the optimal solutions to \req{optimization_problem_z} and \req{optimization_problem_modified} are equivalent under the above rearrangements. 
Thus, we obtain 
\begin{align}
f(k) = \hat{f}(k),\ \forall &f(k) \in {\cal F}_{{\cal I}_k} \ \ \Longleftrightarrow \notag \\ 
 & f_p (k) = \hat{f}_p (k), \ \forall f_p (k) \in \widetilde{\cal F}_{{\cal I}_k}. \label{fkequivalent}
\end{align}
Based on the above notations, the proof of \rlem{main_result_lem} is provided as follows. 

\noindent
\textit{\textbf{(Proof of \rlem{main_result_lem}):}} 
\textbf{(a) $\Rightarrow$ (b) :} 
Suppose that $|{\cal I}_k| < M/2$ holds. 
For the active mode case $a_1(k) = 1$, we have $C_p = C_{p1}$ and since the graph is weakly connected, it is shown from \req{cp2}
that ${\rm rank} (C_p) = nM$ and ${\rm ker} (C_p) = 0$. 
This means that $z_p$ satisfying the constraint $C_p (f_p (k) -z_p) = 0$ is a \textit{single point} and is given by $z_p = f_p (k)$, which shows that the optimization problem in \req{optimization_problem_modified} yields 
$\hat{f}_p (k) = f_p (k)$. This implies that we have $\hat{f} (k) = f (k)$, $\forall f(k) \in {\cal F}_{{\cal I}_k}$, $\forall u(k-1) \in \mathbb{R}^{mM}$. Moreover, since $\hat{f}(k) = f(k)$, it follows that 
\begin{align}
\hat{x}(k) &= A \hat{x}(k-1) + B u(k-1) + \hat{f}(k) \notag \\
             & = A {x}(k-1) + B u(k-1) + {f}(k) \notag \\
             & = x(k), \label{state_correct}
\end{align}
which means that the state $x(k)$ is also correctly estimated. Thus, for the active mode case, the state and fault signals are correctly estimated by solving \req{optimization_problem}. 
Let us now consider the non-active mode case $C_p = C_{p0}$. Let $\widetilde{T}_i = \{i, i+M, i+ 2M, \ldots, i+ (n-1)M\}$, $i\in\{1, \ldots, M\}$ and $\widetilde{T} = \cup_{i\in{\cal I}_k}\widetilde{T}_i$. 
Since $C_p = C_{p0}$ is given by \req{cp} and the graph is weakly connected, we obtain 
\begin{equation}
{\rm ker} (C_p) = [\gamma_1 {1}^\mathsf{T}_M, \gamma_2 {1}^\mathsf{T}_M, \ldots, \gamma_n {1}^\mathsf{T}_M ]^\mathsf{T}, 
\end{equation}
where $\gamma_i \in \mathbb{R}$, $i \in \{1,\ldots,n\}$ and we have used the fact that ${\rm ker}(D^\mathsf{T}) = \gamma {1}^\mathsf{T}_M$ with $\gamma \in\mathbb{R}$ (see \rsec{pre_graphsec}). 
Thus, it follows that 
\begin{equation}
\|v_{\widetilde{T}_i}\|_1 = \sum_{j=1} ^n |\gamma_j|, 
\end{equation}
for any $v \in {\rm ker} (C_p)$. 
Thus, we obtain 
\begin{align}\label{relation0}
\|v_{\widetilde{T}}\|_1 &= \sum_{i\in{\cal I}_k} \|v_{\widetilde{T}_i}\|_1 = |{\cal I}_k| \left (\sum_{j=1} ^n |\gamma_j|\right )
\end{align}
\begin{align}
\|v_{\widetilde{T}^c}\|_1 &= (M - |{\cal I}_k|) \left (\sum_{j=1} ^n |\gamma_j|\right ),  \label{relation1}
\end{align}
where $\widetilde{T}^c = \{1, 2, \ldots, nM\}\backslash \widetilde{T}$. 
Since $|{\cal I}_k| < M/2$, it follows that $\|v_{\widetilde{T}}\|_1 < \|v_{\widetilde{T}^c}\|_1$, $\forall v \in {\rm ker} (C_p) \backslash \{0\}$, i.e., $C_p$ satisfies $\widetilde{T}$-NSP. 
Since $f_p (k) \in \widetilde{\cal F}_{{\cal I}_k}$ implies that $f_p (k)$ is $\widetilde{T}$-sparse, it follows from \rthm{reconstruct} that $f_p (k) = \hat{f}_p(k)$, $\forall f_p (k)\in \widetilde{\cal F}_{{\cal I}_k}$. 
This implies that $\hat{f} (k) = f (k)$, $\forall f(k) \in {\cal F}_{{\cal I}_k}$, $\forall u(k-1) \in \mathbb{R}^{mM}$. From \req{state_correct}, it then follows that $\hat{x}(k) = x(k)$. Hence, {(a) $\Rightarrow$ (b)} holds. 

\textbf{(b) $\Rightarrow$ (a) :} Suppose that (b) holds. To prove by contradiction, suppose that $|{\cal I}_k| \geq M/2$ and consider the non-active mode case $a_1(k) = 0$. Since $|{\cal I}_k| \geq M/2$, it follows from \req{relation0} and \req{relation1} that $\|v_{\widetilde{T}}\|_1 \geq \|v_{\widetilde{T}^c}\|_1$, $\forall v \in {\rm ker} (C_p) \backslash \{0\}$, i.e., $C_p$ does \textit{not} satisfy $\widetilde{T}$-NSP.
Thus, it follows from \rthm{reconstruct} that there exists $f_p (k) \in \widetilde{\cal F}_{{\cal I}_k}$ such that $\hat{f}_p (k) \neq f_p (k)$, which means that $\hat{f} (k) \neq f (k)$ and the optimization problem in \req{optimization_problem} does not provide a correct estimation. 
Indeed, this contradicts to the statement in (b) that $\hat{f} (k) = f(k)$, $\forall f(k) \in {\cal F}_{{\cal I}_k}$. Hence, (b) $\Rightarrow$ (a) holds. \qedwhite \\


\noindent
It is now not difficult to show \rthm{main_result_NSP}:\\

\noindent
\textit{\textbf{(Proof of \rthm{main_result_NSP}):}} 
By \ras{initial_assumption}, the leader is the active mode at $k=0$ ($C= C_1$) and thus ${\rm rank} (C_p) ={\rm rank}(C) = nM$ at $k=0$. This implies that the optimization problem in \req{optimization_problem} yields $\hat{x}(0) = x(0)$, since $C = C_1$ has a trivial kernel and the solution of $x$ satisfying $y(0) = Cx$ is uniquely determined. 
To provide the proof, suppose that (a) holds, i.e., $|{\cal I}_k| < M/2$ for all $k\in\mathbb{N}\backslash \{0\}$. 
Since $\hat{x}(0) = x(0)$, it follows from \rlem{main_result_lem} that $\hat{x} (1) = x(1)$, $\hat{f} (1) = f(1)$, $\forall f(1) \in {\mathbb{R}}^{nM} _{{\cal I}_1}$, $\forall u(0) \in \mathbb{R}^{mM}$. Since $\hat{x} (1) = x(1)$, it then follows that $\hat{x} (2)= x(2)$, $\hat{f} (2) = f(2)$, $\forall f(2) \in {\mathbb{R}}^{nM} _{{\cal I}_2}$, $\forall u(1) \in \mathbb{R}^{mM}$. Similarly, it follows recursively by applying \rlem{main_result_lem} that $\hat{x} (k) = x(k)$, $\hat{f} (k) = f(k)$, $\forall f(k) \in {\cal F}_{{\cal I}_k}$, $\forall u(k-1) \in \mathbb{R}^{mM}$. Hence, (a) $\Rightarrow$ (b) holds. 
Conversely, suppose that (b) holds. 
Since $\hat{x}(0) = x(0)$, it follows from \rlem{main_result_lem} that $|{\cal I}_1| < M/2$. Similarly, it follows recursively by applying \rlem{main_result_lem} that $|{\cal I}_k| < M/2$, $\forall k\in\mathbb{N}\backslash \{0\}$. Hence, (b) $\Rightarrow$ (a) holds. \qedwhite

\section{Proof of \rthm{main_result_error}}\label{proof_appendix2}
Similarly to the proof of \rlem{main_result_lem} and \rthm{main_result_NSP}, it is required to take several steps to modify the optmization problem in \req{optimization_problem_simple}. 
Since $\|x (k-1) - \hat{x} (k-1)\|_1 \leq d_{\max} $, $x(k-1)$ is expressed as $x (k-1) = \hat{x} (k-1) + \tilde{d}$, where $\|\tilde{d}\|_1 \leq d_{\max}$. Thus we obtain $y(k) = C(A \hat{x} (k-1) + A \tilde{d} + B u(k-1) + f(k))$ and \req{optimization_problem_simple} becomes 
\begin{equation}\label{optimization_problem_error}
\underset{z}{{\min}}\ \| z \| _1 \ \ \ {\rm s.t.,}\ \ Cz = C (f(k) - d), 
\end{equation}
where we let $z = x - A \hat{x}(k-1) - Bu(k-1)$ and $d = A \tilde{d}$. Thus, \req{optimization_problem_error} differs from \req{optimization_problem_z} in the sense that the term $d$ arised from the estimation error at $k-1$ is added in the constraint. Since $\eta = \sum^{nM} _{i=1} \sum^{nM} _{j=1} |A^{(i, j)}|$, it follows that $\|d\|_1 = \|A\tilde{d}\|_1 \leq \eta d_{\max}$. By rearranging the vectors $z, f(k), d$ and the matrix $C$ as with the procedure presented in Appendix~B, the problem becomes 
\begin{equation}\label{optimization_problem_error_modified}
\underset{z}{{\min}}\ \| z_p \| _1 \ \ \ {\rm s.t.,}\ \ C_p z_p = C_p (f_p (k) - d_p), 
\end{equation}
where $d_p = [d_{p_1}^\mathsf{T}, d_{p_2}^\mathsf{T}, \ldots, d_{p_n}^\mathsf{T}]^\mathsf{T}$ with $d_{p_i} = d_{T_i}$, $i\in\{1, 2, \ldots, n\}$. In \req{optimization_problem_error_modified}, $C_p = C_{p0}$ if $a_1 (k) = 0$ and  $C_p = C_{p1}$ if $a_1 (k) = 1$. It follows that $\|d_p\|_1 \leq \eta d_{\max}$ and 
\begin{align}
\|f(k) - & \hat{f}(k)\|_1 \leq \epsilon,\ \forall f(k) \in {\cal F}_{{\cal I}_k}\Longleftrightarrow \notag \\
& \ \ \|f_p (k) - \hat{f}_p (k)\|_1 \leq \epsilon, \ \forall f_p (k) \in \widetilde{\cal F}_{{\cal I}_k}, 
\end{align}
for a given $\epsilon > 0$. 

Now, let $\hat{f}_p(k)$ be the optimal solution of $z_p$ to \req{optimization_problem_error_modified}. For the active mode case, we have $C_p = C_{p1}$ and it follows that ${\rm rank} (C_p) = nM$ and ${\rm ker} (C_p) = 0$. Thus, we have $\hat{f}_p(k) = f_p(k) - d_p$ and $\|f_p (k) - \hat{f}_p (k) \|_1 \leq \eta d_{\max}$, which means that $\|f (k) - \hat{f} (k) \|_1 \leq \eta d_{\max}$. Since $\eta d_{\max} < \frac{2(M-|{\cal I}_k|)}{M-2|{\cal I}_k|} \eta d _{\max}$, we obtain \req{error_bound}. 
Consider now the nonactive mode case, i.e., $C_p = C_{p0}$. 
Since $\hat{f}_p (k)$ is the optimal solution, it follows that $\|\hat{f}_p (k)\| \leq \|f_p(k)\|$. Since $C_p \hat{f}_p (k) = C_p (f_p(k) - d_p)$, there exists $v\in {\rm ker} (C_p)$ such that $\hat{f}_p (k) = f_p (k) - d_p + v$. Thus, it follows that 
\begin{align}
&\|f_p \|_1 \geq \| f_p  - d_p + v\|_1 \notag \\
               & = \sum_{i\in \widetilde{T}} |f^{(i)}_p - d^{(i)}_p + v^{(i)}| + \sum_{i\in \widetilde{T}^c} |f^{(i)}_p - d^{(i)}_p + v^{(i)}| \notag \\ 
               & \geq \sum_{i\in \widetilde{T}} \left (|f^{(i)}_p - d^{(i)}_p| - |v^{(i)}|\right ) + \sum_{i\in \widetilde{T}^c} \left (|v^{(i)}| - |f^{(i)}_p - d^{(i)}_p|\right )  \notag 
\end{align}
where the time index $k$ for $f_p$ is omitted for brevity, and $\widetilde{T} = \cup_{i\in{\cal I}_k}\widetilde{T}_i$ with $\widetilde{T}_i = \{i, i+M, i+ 2M, \ldots, i+ (n-1)M\}$, $i\in\{1, \ldots, M\}$. For the third inequality in the above, we used the triangle inequality ($|a + b| \geq |a| - |b|$). 
Thus, it follows that 
\begin{align}
\|f_p \|_1 &\geq  \| (f_p)_{\widetilde{T}} - (d_p)_{\widetilde{T}}\|_1  - \|(d_p)_{\widetilde{T}}\|_1 - \|v_{\widetilde{T}}\|_1 + \|v_{\widetilde{T}^c}\|_1, \notag 
\end{align}
where we used $f^{(i)}_p={0}$, $\forall i\in \widetilde{T}^c$ as $f_p\in \widetilde{\cal F}_{{\cal I}_k}$. 
Moreover, it follows from \req{relation1} that 
$\|v_{\widetilde{T}}\|_1 = \frac{|{\cal I}_k|}{M-|{\cal I}_k|}\ \|v_{\widetilde{T}^c}\|_1$. 
By plugging this into the above inequality, we obtain 
\begin{align}
\|v_{\widetilde{T}^c}\|_1 \leq \frac{M-|{\cal I}_k|}{M-2 |{\cal I}_k|} ( \|f_p\|_1 & - \|(f_p)_{\widetilde{T}} - (d_p)_{\widetilde{T}}\|_1 \notag \\ 
& + \|(d_p)_{\widetilde{T}^c}\|_1). 
\end{align} 
Since we have $\|f_p\|_1 = \|(f_p)_{\widetilde{T}}\| + \|(f_p)_{\widetilde{T}^c}\| =\|(f_p)_{\widetilde{T}}\|$ and by using the triangle inequality, we obtain 
\begin{align}
\|v_{\widetilde{T}^c}\|_1 &\leq \cfrac{M-|{\cal I}_k|}{M-2 |{\cal I}_k|} \left( \|(d_p)_{\tilde{T}}\|_1 + \|(d_p)_{\tilde{T}^c}\|_1\right) \notag \\
&\leq \cfrac{M-|{\cal I}_k|}{M-2 |{\cal I}_k|}\ \eta d_{\max}. \notag 
\end{align} 
Therefore, it follows that 
\begin{align}
\|\hat{f}_p (k) - f_p (k)\|_1 &= \|-d_p + v\|_1 \leq \eta d_{\max} + \|v\|_1 \notag \\ 
                              &= \eta d_{\max} + \left ( 1 + \cfrac{|{\cal I}_k|} {M- |{\cal I}_k|} \right ) \|v_{\widetilde{T}^c}\|_1 \notag \\
                              &\leq \cfrac{2(M-|{\cal I}_k|)} {M- 2|{\cal I}_k|}\ \eta d_{\max}, 
\end{align} 
which implies that \req{error_bound} holds.  \qedwhite \\ 

\end{document}